\documentclass[12pt]{article}

\usepackage[margin=1in]{geometry}
\usepackage{amsmath,amssymb,amsfonts,amsthm}
\usepackage{graphicx}
\usepackage{xcolor,tikz}
\usepackage{hyperref}
\usetikzlibrary{positioning,arrows.meta,calc,matrix,decorations}

\numberwithin{equation}{section}

\newtheorem{theorem}[equation]{Theorem}
\newtheorem{proposition}[equation]{Proposition}
\newtheorem{claim}[equation]{Claim}

\theoremstyle{remark}
\newtheorem{remark}[equation]{Remark}

\title{On Alternating 6-Cycles in Edge-Coloured Graphs}

\author{Hao Chen\thanks{School of Mathematical Sciences, Soochow University,
Suzhou 215006, China. E-mail:
\href{mailto:chenhao@suda.edu.cn}{\texttt{chenhao@suda.edu.cn}}.}
\and
Jonathan A. Noel\thanks{Department of Mathematics and Statistics,
University of Victoria, Victoria, B.C., Canada. E-mail:
\href{mailto:noelj@uvic.ca}{\texttt{noelj@uvic.ca}}.}}

\date{}

\newcommand\scalemath[2]{%
  \scalebox{#1}{\mbox{\ensuremath{\displaystyle #2}}}%
}

\DeclareTextCompositeCommand{\v}{OT1}{l}{l\nobreak\hspace{-.1em}'}
\DeclareTextCompositeCommand{\v}{OT1}{t}{t\nobreak\hspace{-.1em}'%
  \nobreak\hspace{-.15em}}

\DeclareMathOperator{\inj}{inj}
\DeclareMathOperator{\bip}{bip}

\DeclareMathOperator{\aut}{aut}

\begin{document}

\maketitle

\begin{abstract}
In this short note, we use flag algebras to prove that the number of colour-alternating $6$-cycles in a red/blue colouring of a large clique is asymptotically maximized by a uniformly random colouring. This settles the first open case of a problem of Basit, Granet, Horsley, K\"undgen and Staden. 
\end{abstract}

\section{Introduction}

An \emph{edge-coloured graph} is a graph whose edges are coloured red or blue. For $k\geq2$, the \emph{alternating cycle} $C_{2k}^A$ is the edge-coloured graph consisting of the cycle graph $C_{2k}$ where every vertex is incident to exactly one edge of each colour. 

Our goal is to determine the maximum number of ``copies'' of the alternating $6$-cycle in a large edge-coloured graph $G$. Here, when speaking of ``copies,'' we mean ``homomorphic copies.'' More precisely, a \emph{homomorphism} from an edge-coloured graph $H$ to an edge-coloured graph $G$ is a map $f:V(H)\to V(G)$ with the property that if $uv$ is an edge of $H$, then $f(u)f(v)$ is an edge of $G$ with the same colour as $uv$. Let $\hom(H,G)$ be the number of homomorphisms from $H$ to $G$. The \emph{homomorphism density} of $H$ in $G$ is $t(H,G)=\hom(H,G)/v(G)^{v(H)}$. The following theorem solves the first case of a problem posed recently by Basit et al.~\cite[Problem~9.1]{Basit+25+}.

\begin{theorem}
\label{th:main}
Every edge-coloured graph $G$ satisfies $t(C_6^A,G)\leq (1/2)^6$.
\end{theorem}

It is not hard to see that this is asymptotically tight by taking $G$ to be a uniformly random edge-coloured clique on a large number of vertices. Our proof of Theorem~\ref{th:main}, provided in the next section, is a relatively simple application of the flag algebra method of Razborov~\cite{Razborov07}.

The problem of determining the maximum density of copies of an edge-coloured graph $H$ in an edge-coloured graph $G$, called the \emph{semi-inducibility problem}, was introduced recently by Basit et al.~\cite{Basit+25+}, who solved it for several classes of graphs including alternating paths of even length and alternating cycles of length divisible by four. In a recent paper together with Clemen, we used entropy to settle the case of alternating paths of odd length~\cite{ChenClemenNoel25+}. It would be interesting to solve the semi-inducibility problem for alternating cycles of length $2\bmod 4$ in full generality~\cite[Problem~9.1]{Basit+25+}. The method in the present paper is, unfortunately, unlikely to lead to a general proof, but it does suggest that the correct bound may be of the form $t(C_{4k+2}^A,G)\leq (1/2)^{4k+2}$. 

\begin{remark}
We are aware that another group of researchers, Balogh, Lidick\'y, Mubayi and Pfender (personal communication), independently proved Theorem~\ref{th:main} using flag algebras and have obtained other related results. 
\end{remark}

\section{The Proof}

Our first step in the proof of Theorem~\ref{th:main} is to reduce it to a statement about injective homomorphisms. For edge-coloured graphs $H$ and $G$ with $v(G)\geq v(H)$, let $\hom_{\inj}(H,G)$ be the number of injective homomorphisms from $H$ to $G$; when $v(G)<v(H)$, then we let $\hom_{\inj}(H,G):=0$. Also, define $t_{\inj}(H,G):=\frac{(v(G)-v(H))!}{v(G)!}\hom_{\inj}(H,G)$. That is, $t_{\inj}(H,G)$ is the probability that a uniformly random injective function from $V(H)$ to $V(G)$ is a homomorphism. We show that Theorem~\ref{th:main} can be deduced from the following theorem; after this, we focus the rest of the paper on proving Theorem~\ref{th:injective}. 

\begin{theorem}
\label{th:injective}
Every edge-coloured graph $G$ satisfies $t_{\inj}(C_6^A,G)\leq (1/2)^6 + o(1)$.
\end{theorem}

\begin{proof}[Proof of Theorem~\ref{th:main}, assuming Theorem~\ref{th:injective}]
Suppose that there exists an edge-coloured graph $G$ such that $t(C_6^A,G)>(1/2)^6$. For each $N\geq1$, let $G_N$ be the edge-coloured graph obtained by ``blowing up'' each vertex of $G$ into a set of size $N$. That is, $G_N$ is the edge-coloured graph consisting of $N\cdot v(G)$ vertices divided into $v(G)$ classes of size $N$, one for each vertex of $G$, where the edges between two classes are the same colour as the edge between the corresponding vertices of $G$. It is easily observed that $\lim_{N\to\infty} t_{\inj}(C_6^A,G_N)=t(C_6^A,G)$ and so, for large enough $N$, the graph $G_N$ contradicts Theorem~\ref{th:injective}. 
\end{proof}

We prove Theorem~\ref{th:injective}. Let $G$ be an edge-coloured graph on $n$ vertices. Our aim is to prove that $t_{\inj}(C_6^A,G)\leq (1/2)^6+o(1)$. We may assume that $G$ is an edge-coloured clique; this is because adding any missing edge with an arbitrary colour cannot decrease $t_{\inj}(C_6^A,G)$. 

There are $2^9=512$ distinct red/blue edge colourings of the complete bipartite graph $K_{3,3}$. Let $H_1,\dots,H_{512}$ be the corresponding edge-coloured graphs where, for convenience, we assume $V(K_{3,3})=V(H_1)=\cdots=V(H_{512})$. Since $G$ is a clique, for any injective function $f:V(K_{3,3})\to V(G)$, there is a unique index $i\in\{1,\dots,512\}$ such that $f$ is a homomorphism of $H_i$ to $G$. Thus,
\begin{equation}
\label{eq:sumToOne}
    \sum_{i=1}^{512}t_{\inj}(H_i,G) = 1.
\end{equation}

\begin{figure}[htbp]
    \centering
\begin{tikzpicture}[scale=0.3]
\node[circle, draw, fill=black, inner sep=1.3pt] (A1) at (3,0) {};
\node[circle, draw, fill=black, inner sep=1.3pt] (A2) at (3,2) {};
\node[circle, draw, fill=black, inner sep=1.3pt] (A3) at (3,4) {};
\node[circle, draw, fill=black, inner sep=1.3pt] (B1) at (0,0) {};
\node[circle, draw, fill=black, inner sep=1.3pt] (B2) at (0,2) {};
\node[circle, draw, fill=black, inner sep=1.3pt] (B3) at (0,4) {};
\node (J) at (1.5,-1) {$J_{1}$};
\draw[red, thick] (A1) -- (B1);
\draw[red, thick] (A1) -- (B2);
\draw[red, thick] (A1) -- (B3);
\draw[red, thick] (A2) -- (B1);
\draw[red, thick] (A2) -- (B2);
\draw[red, thick] (A2) -- (B3);
\draw[red, thick] (A3) -- (B1);
\draw[red, thick] (A3) -- (B2);
\draw[red, thick] (A3) -- (B3);
\end{tikzpicture}
\hspace{1em}
\begin{tikzpicture}[scale=0.3]
\node[circle, draw, fill=black, inner sep=1.3pt] (A1) at (3,0) {};
\node[circle, draw, fill=black, inner sep=1.3pt] (A2) at (3,2) {};
\node[circle, draw, fill=black, inner sep=1.3pt] (A3) at (3,4) {};
\node[circle, draw, fill=black, inner sep=1.3pt] (B1) at (0,0) {};
\node[circle, draw, fill=black, inner sep=1.3pt] (B2) at (0,2) {};
\node[circle, draw, fill=black, inner sep=1.3pt] (B3) at (0,4) {};
\node (J) at (1.5,-1) {$J_{2}$};
\draw[blue, dashed, thick] (A1) -- (B1);
\draw[red, thick] (A1) -- (B2);
\draw[red, thick] (A1) -- (B3);
\draw[red, thick] (A2) -- (B1);
\draw[red, thick] (A2) -- (B2);
\draw[red, thick] (A2) -- (B3);
\draw[red, thick] (A3) -- (B1);
\draw[red, thick] (A3) -- (B2);
\draw[red, thick] (A3) -- (B3);
\end{tikzpicture}
\hspace{1em}
\begin{tikzpicture}[scale=0.3]
\node[circle, draw, fill=black, inner sep=1.3pt] (A1) at (3,0) {};
\node[circle, draw, fill=black, inner sep=1.3pt] (A2) at (3,2) {};
\node[circle, draw, fill=black, inner sep=1.3pt] (A3) at (3,4) {};
\node[circle, draw, fill=black, inner sep=1.3pt] (B1) at (0,0) {};
\node[circle, draw, fill=black, inner sep=1.3pt] (B2) at (0,2) {};
\node[circle, draw, fill=black, inner sep=1.3pt] (B3) at (0,4) {};
\node (J) at (1.5,-1) {$J_{3}$};
\draw[blue, dashed, thick] (A1) -- (B1);
\draw[red, thick] (A1) -- (B2);
\draw[red, thick] (A1) -- (B3);
\draw[red, thick] (A2) -- (B1);
\draw[blue, dashed, thick] (A2) -- (B2);
\draw[red, thick] (A2) -- (B3);
\draw[red, thick] (A3) -- (B1);
\draw[red, thick] (A3) -- (B2);
\draw[red, thick] (A3) -- (B3);
\end{tikzpicture}
\hspace{1em}
\begin{tikzpicture}[scale=0.3]
\node[circle, draw, fill=black, inner sep=1.3pt] (A1) at (3,0) {};
\node[circle, draw, fill=black, inner sep=1.3pt] (A2) at (3,2) {};
\node[circle, draw, fill=black, inner sep=1.3pt] (A3) at (3,4) {};
\node[circle, draw, fill=black, inner sep=1.3pt] (B1) at (0,0) {};
\node[circle, draw, fill=black, inner sep=1.3pt] (B2) at (0,2) {};
\node[circle, draw, fill=black, inner sep=1.3pt] (B3) at (0,4) {};
\node (J) at (1.5,-1) {$J_{4}$};
\draw[blue, dashed, thick] (A1) -- (B1);
\draw[red, thick] (A1) -- (B2);
\draw[red, thick] (A1) -- (B3);
\draw[red, thick] (A2) -- (B1);
\draw[blue, dashed, thick] (A2) -- (B2);
\draw[red, thick] (A2) -- (B3);
\draw[red, thick] (A3) -- (B1);
\draw[red, thick] (A3) -- (B2);
\draw[blue, dashed, thick] (A3) -- (B3);
\end{tikzpicture}
\hspace{1em}
\begin{tikzpicture}[scale=0.3]
\node[circle, draw, fill=black, inner sep=1.3pt] (A1) at (3,0) {};
\node[circle, draw, fill=black, inner sep=1.3pt] (A2) at (3,2) {};
\node[circle, draw, fill=black, inner sep=1.3pt] (A3) at (3,4) {};
\node[circle, draw, fill=black, inner sep=1.3pt] (B1) at (0,0) {};
\node[circle, draw, fill=black, inner sep=1.3pt] (B2) at (0,2) {};
\node[circle, draw, fill=black, inner sep=1.3pt] (B3) at (0,4) {};
\node (J) at (1.5,-1) {$J_{5}$};
\draw[blue, dashed, thick] (A1) -- (B1);
\draw[blue, dashed, thick] (A1) -- (B2);
\draw[red, thick] (A1) -- (B3);
\draw[red, thick] (A2) -- (B1);
\draw[red, thick] (A2) -- (B2);
\draw[red, thick] (A2) -- (B3);
\draw[red, thick] (A3) -- (B1);
\draw[red, thick] (A3) -- (B2);
\draw[red, thick] (A3) -- (B3);
\end{tikzpicture}
\hspace{1em}
\begin{tikzpicture}[scale=0.3]
\node[circle, draw, fill=black, inner sep=1.3pt] (A1) at (3,0) {};
\node[circle, draw, fill=black, inner sep=1.3pt] (A2) at (3,2) {};
\node[circle, draw, fill=black, inner sep=1.3pt] (A3) at (3,4) {};
\node[circle, draw, fill=black, inner sep=1.3pt] (B1) at (0,0) {};
\node[circle, draw, fill=black, inner sep=1.3pt] (B2) at (0,2) {};
\node[circle, draw, fill=black, inner sep=1.3pt] (B3) at (0,4) {};
\node (J) at (1.5,-1) {$J_{6}$};
\draw[blue, dashed, thick] (A1) -- (B1);
\draw[blue, dashed, thick] (A1) -- (B2);
\draw[red, thick] (A1) -- (B3);
\draw[red, thick] (A2) -- (B1);
\draw[red, thick] (A2) -- (B2);
\draw[blue, dashed, thick] (A2) -- (B3);
\draw[red, thick] (A3) -- (B1);
\draw[red, thick] (A3) -- (B2);
\draw[red, thick] (A3) -- (B3);
\end{tikzpicture}
\hspace{1em}
\begin{tikzpicture}[scale=0.3]
\node[circle, draw, fill=black, inner sep=1.3pt] (A1) at (3,0) {};
\node[circle, draw, fill=black, inner sep=1.3pt] (A2) at (3,2) {};
\node[circle, draw, fill=black, inner sep=1.3pt] (A3) at (3,4) {};
\node[circle, draw, fill=black, inner sep=1.3pt] (B1) at (0,0) {};
\node[circle, draw, fill=black, inner sep=1.3pt] (B2) at (0,2) {};
\node[circle, draw, fill=black, inner sep=1.3pt] (B3) at (0,4) {};
\node (J) at (1.5,-1) {$J_{7}$};
\draw[blue, dashed, thick] (A1) -- (B1);
\draw[blue, dashed, thick] (A1) -- (B2);
\draw[red, thick] (A1) -- (B3);
\draw[red, thick] (A2) -- (B1);
\draw[red, thick] (A2) -- (B2);
\draw[blue, dashed, thick] (A2) -- (B3);
\draw[red, thick] (A3) -- (B1);
\draw[red, thick] (A3) -- (B2);
\draw[blue, dashed, thick] (A3) -- (B3);
\end{tikzpicture}
\hspace{1em}
\begin{tikzpicture}[scale=0.3]
\node[circle, draw, fill=black, inner sep=1.3pt] (A1) at (3,0) {};
\node[circle, draw, fill=black, inner sep=1.3pt] (A2) at (3,2) {};
\node[circle, draw, fill=black, inner sep=1.3pt] (A3) at (3,4) {};
\node[circle, draw, fill=black, inner sep=1.3pt] (B1) at (0,0) {};
\node[circle, draw, fill=black, inner sep=1.3pt] (B2) at (0,2) {};
\node[circle, draw, fill=black, inner sep=1.3pt] (B3) at (0,4) {};
\node (J) at (1.5,-1) {$J_{8}$};
\draw[blue, dashed, thick] (A1) -- (B1);
\draw[blue, dashed, thick] (A1) -- (B2);
\draw[red, thick] (A1) -- (B3);
\draw[blue, dashed, thick] (A2) -- (B1);
\draw[red, thick] (A2) -- (B2);
\draw[red, thick] (A2) -- (B3);
\draw[red, thick] (A3) -- (B1);
\draw[red, thick] (A3) -- (B2);
\draw[red, thick] (A3) -- (B3);
\end{tikzpicture}
\hspace{1em}
\begin{tikzpicture}[scale=0.3]
\node[circle, draw, fill=black, inner sep=1.3pt] (A1) at (3,0) {};
\node[circle, draw, fill=black, inner sep=1.3pt] (A2) at (3,2) {};
\node[circle, draw, fill=black, inner sep=1.3pt] (A3) at (3,4) {};
\node[circle, draw, fill=black, inner sep=1.3pt] (B1) at (0,0) {};
\node[circle, draw, fill=black, inner sep=1.3pt] (B2) at (0,2) {};
\node[circle, draw, fill=black, inner sep=1.3pt] (B3) at (0,4) {};
\node (J) at (1.5,-1) {$J_{9}$};
\draw[blue, dashed, thick] (A1) -- (B1);
\draw[blue, dashed, thick] (A1) -- (B2);
\draw[red, thick] (A1) -- (B3);
\draw[blue, dashed, thick] (A2) -- (B1);
\draw[red, thick] (A2) -- (B2);
\draw[red, thick] (A2) -- (B3);
\draw[red, thick] (A3) -- (B1);
\draw[red, thick] (A3) -- (B2);
\draw[blue, dashed, thick] (A3) -- (B3);
\end{tikzpicture}

\begin{tikzpicture}[scale=0.3]
\node[circle, draw, fill=black, inner sep=1.3pt] (A1) at (3,0) {};
\node[circle, draw, fill=black, inner sep=1.3pt] (A2) at (3,2) {};
\node[circle, draw, fill=black, inner sep=1.3pt] (A3) at (3,4) {};
\node[circle, draw, fill=black, inner sep=1.3pt] (B1) at (0,0) {};
\node[circle, draw, fill=black, inner sep=1.3pt] (B2) at (0,2) {};
\node[circle, draw, fill=black, inner sep=1.3pt] (B3) at (0,4) {};
\node (J) at (1.5,-1) {$J_{10}$};
\draw[blue, dashed, thick] (A1) -- (B1);
\draw[blue, dashed, thick] (A1) -- (B2);
\draw[red, thick] (A1) -- (B3);
\draw[blue, dashed, thick] (A2) -- (B1);
\draw[red, thick] (A2) -- (B2);
\draw[blue, dashed, thick] (A2) -- (B3);
\draw[red, thick] (A3) -- (B1);
\draw[red, thick] (A3) -- (B2);
\draw[red, thick] (A3) -- (B3);
\end{tikzpicture}
\hspace{1em}
\begin{tikzpicture}[scale=0.3]
\node[circle, draw, fill=black, inner sep=1.3pt] (A1) at (3,0) {};
\node[circle, draw, fill=black, inner sep=1.3pt] (A2) at (3,2) {};
\node[circle, draw, fill=black, inner sep=1.3pt] (A3) at (3,4) {};
\node[circle, draw, fill=black, inner sep=1.3pt] (B1) at (0,0) {};
\node[circle, draw, fill=black, inner sep=1.3pt] (B2) at (0,2) {};
\node[circle, draw, fill=black, inner sep=1.3pt] (B3) at (0,4) {};
\node (J) at (1.5,-1) {$J_{11}$};
\draw[blue, dashed, thick] (A1) -- (B1);
\draw[blue, dashed, thick] (A1) -- (B2);
\draw[red, thick] (A1) -- (B3);
\draw[blue, dashed, thick] (A2) -- (B1);
\draw[red, thick] (A2) -- (B2);
\draw[blue, dashed, thick] (A2) -- (B3);
\draw[red, thick] (A3) -- (B1);
\draw[blue, dashed, thick] (A3) -- (B2);
\draw[red, thick] (A3) -- (B3);
\end{tikzpicture}
\hspace{1em}
\begin{tikzpicture}[scale=0.3]
\node[circle, draw, fill=black, inner sep=1.3pt] (A1) at (3,0) {};
\node[circle, draw, fill=black, inner sep=1.3pt] (A2) at (3,2) {};
\node[circle, draw, fill=black, inner sep=1.3pt] (A3) at (3,4) {};
\node[circle, draw, fill=black, inner sep=1.3pt] (B1) at (0,0) {};
\node[circle, draw, fill=black, inner sep=1.3pt] (B2) at (0,2) {};
\node[circle, draw, fill=black, inner sep=1.3pt] (B3) at (0,4) {};
\node (J) at (1.5,-1) {$J_{12}$};
\draw[blue, dashed, thick] (A1) -- (B1);
\draw[blue, dashed, thick] (A1) -- (B2);
\draw[red, thick] (A1) -- (B3);
\draw[blue, dashed, thick] (A2) -- (B1);
\draw[red, thick] (A2) -- (B2);
\draw[blue, dashed, thick] (A2) -- (B3);
\draw[red, thick] (A3) -- (B1);
\draw[blue, dashed, thick] (A3) -- (B2);
\draw[blue, dashed, thick] (A3) -- (B3);
\end{tikzpicture}
\hspace{1em}
\begin{tikzpicture}[scale=0.3]
\node[circle, draw, fill=black, inner sep=1.3pt] (A1) at (3,0) {};
\node[circle, draw, fill=black, inner sep=1.3pt] (A2) at (3,2) {};
\node[circle, draw, fill=black, inner sep=1.3pt] (A3) at (3,4) {};
\node[circle, draw, fill=black, inner sep=1.3pt] (B1) at (0,0) {};
\node[circle, draw, fill=black, inner sep=1.3pt] (B2) at (0,2) {};
\node[circle, draw, fill=black, inner sep=1.3pt] (B3) at (0,4) {};
\node (J) at (1.5,-1) {$J_{13}$};
\draw[blue, dashed, thick] (A1) -- (B1);
\draw[blue, dashed, thick] (A1) -- (B2);
\draw[red, thick] (A1) -- (B3);
\draw[blue, dashed, thick] (A2) -- (B1);
\draw[blue, dashed, thick] (A2) -- (B2);
\draw[red, thick] (A2) -- (B3);
\draw[red, thick] (A3) -- (B1);
\draw[red, thick] (A3) -- (B2);
\draw[red, thick] (A3) -- (B3);
\end{tikzpicture}
\hspace{1em}
\begin{tikzpicture}[scale=0.3]
\node[circle, draw, fill=black, inner sep=1.3pt] (A1) at (3,0) {};
\node[circle, draw, fill=black, inner sep=1.3pt] (A2) at (3,2) {};
\node[circle, draw, fill=black, inner sep=1.3pt] (A3) at (3,4) {};
\node[circle, draw, fill=black, inner sep=1.3pt] (B1) at (0,0) {};
\node[circle, draw, fill=black, inner sep=1.3pt] (B2) at (0,2) {};
\node[circle, draw, fill=black, inner sep=1.3pt] (B3) at (0,4) {};
\node (J) at (1.5,-1) {$J_{14}$};
\draw[blue, dashed, thick] (A1) -- (B1);
\draw[blue, dashed, thick] (A1) -- (B2);
\draw[red, thick] (A1) -- (B3);
\draw[blue, dashed, thick] (A2) -- (B1);
\draw[blue, dashed, thick] (A2) -- (B2);
\draw[red, thick] (A2) -- (B3);
\draw[red, thick] (A3) -- (B1);
\draw[red, thick] (A3) -- (B2);
\draw[blue, dashed, thick] (A3) -- (B3);
\end{tikzpicture}
\hspace{1em}
\begin{tikzpicture}[scale=0.3]
\node[circle, draw, fill=black, inner sep=1.3pt] (A1) at (3,0) {};
\node[circle, draw, fill=black, inner sep=1.3pt] (A2) at (3,2) {};
\node[circle, draw, fill=black, inner sep=1.3pt] (A3) at (3,4) {};
\node[circle, draw, fill=black, inner sep=1.3pt] (B1) at (0,0) {};
\node[circle, draw, fill=black, inner sep=1.3pt] (B2) at (0,2) {};
\node[circle, draw, fill=black, inner sep=1.3pt] (B3) at (0,4) {};
\node (J) at (1.5,-1) {$J_{15}$};
\draw[blue, dashed, thick] (A1) -- (B1);
\draw[blue, dashed, thick] (A1) -- (B2);
\draw[blue, dashed, thick] (A1) -- (B3);
\draw[red, thick] (A2) -- (B1);
\draw[red, thick] (A2) -- (B2);
\draw[red, thick] (A2) -- (B3);
\draw[red, thick] (A3) -- (B1);
\draw[red, thick] (A3) -- (B2);
\draw[red, thick] (A3) -- (B3);
\end{tikzpicture}
\hspace{1em}
\begin{tikzpicture}[scale=0.3]
\node[circle, draw, fill=black, inner sep=1.3pt] (A1) at (3,0) {};
\node[circle, draw, fill=black, inner sep=1.3pt] (A2) at (3,2) {};
\node[circle, draw, fill=black, inner sep=1.3pt] (A3) at (3,4) {};
\node[circle, draw, fill=black, inner sep=1.3pt] (B1) at (0,0) {};
\node[circle, draw, fill=black, inner sep=1.3pt] (B2) at (0,2) {};
\node[circle, draw, fill=black, inner sep=1.3pt] (B3) at (0,4) {};
\node (J) at (1.5,-1) {$J_{16}$};
\draw[blue, dashed, thick] (A1) -- (B1);
\draw[blue, dashed, thick] (A1) -- (B2);
\draw[blue, dashed, thick] (A1) -- (B3);
\draw[blue, dashed, thick] (A2) -- (B1);
\draw[red, thick] (A2) -- (B2);
\draw[red, thick] (A2) -- (B3);
\draw[red, thick] (A3) -- (B1);
\draw[red, thick] (A3) -- (B2);
\draw[red, thick] (A3) -- (B3);
\end{tikzpicture}
\hspace{1em}
\begin{tikzpicture}[scale=0.3]
\node[circle, draw, fill=black, inner sep=1.3pt] (A1) at (3,0) {};
\node[circle, draw, fill=black, inner sep=1.3pt] (A2) at (3,2) {};
\node[circle, draw, fill=black, inner sep=1.3pt] (A3) at (3,4) {};
\node[circle, draw, fill=black, inner sep=1.3pt] (B1) at (0,0) {};
\node[circle, draw, fill=black, inner sep=1.3pt] (B2) at (0,2) {};
\node[circle, draw, fill=black, inner sep=1.3pt] (B3) at (0,4) {};
\node (J) at (1.5,-1) {$J_{17}$};
\draw[blue, dashed, thick] (A1) -- (B1);
\draw[blue, dashed, thick] (A1) -- (B2);
\draw[blue, dashed, thick] (A1) -- (B3);
\draw[blue, dashed, thick] (A2) -- (B1);
\draw[red, thick] (A2) -- (B2);
\draw[red, thick] (A2) -- (B3);
\draw[red, thick] (A3) -- (B1);
\draw[blue, dashed, thick] (A3) -- (B2);
\draw[red, thick] (A3) -- (B3);
\end{tikzpicture}
\hspace{1em}
\begin{tikzpicture}[scale=0.3]
\node[circle, draw, fill=black, inner sep=1.3pt] (A1) at (3,0) {};
\node[circle, draw, fill=black, inner sep=1.3pt] (A2) at (3,2) {};
\node[circle, draw, fill=black, inner sep=1.3pt] (A3) at (3,4) {};
\node[circle, draw, fill=black, inner sep=1.3pt] (B1) at (0,0) {};
\node[circle, draw, fill=black, inner sep=1.3pt] (B2) at (0,2) {};
\node[circle, draw, fill=black, inner sep=1.3pt] (B3) at (0,4) {};
\node (J) at (1.5,-1) {$J_{18}$};
\draw[blue, dashed, thick] (A1) -- (B1);
\draw[blue, dashed, thick] (A1) -- (B2);
\draw[blue, dashed, thick] (A1) -- (B3);
\draw[blue, dashed, thick] (A2) -- (B1);
\draw[red, thick] (A2) -- (B2);
\draw[red, thick] (A2) -- (B3);
\draw[blue, dashed, thick] (A3) -- (B1);
\draw[red, thick] (A3) -- (B2);
\draw[red, thick] (A3) -- (B3);
\end{tikzpicture}

\begin{tikzpicture}[scale=0.3]
\node[circle, draw, fill=black, inner sep=1.3pt] (A1) at (3,0) {};
\node[circle, draw, fill=black, inner sep=1.3pt] (A2) at (3,2) {};
\node[circle, draw, fill=black, inner sep=1.3pt] (A3) at (3,4) {};
\node[circle, draw, fill=black, inner sep=1.3pt] (B1) at (0,0) {};
\node[circle, draw, fill=black, inner sep=1.3pt] (B2) at (0,2) {};
\node[circle, draw, fill=black, inner sep=1.3pt] (B3) at (0,4) {};
\node (J) at (1.5,-1) {$J_{19}$};
\draw[blue, dashed, thick] (A1) -- (B1);
\draw[blue, dashed, thick] (A1) -- (B2);
\draw[blue, dashed, thick] (A1) -- (B3);
\draw[blue, dashed, thick] (A2) -- (B1);
\draw[blue, dashed, thick] (A2) -- (B2);
\draw[red, thick] (A2) -- (B3);
\draw[red, thick] (A3) -- (B1);
\draw[red, thick] (A3) -- (B2);
\draw[red, thick] (A3) -- (B3);
\end{tikzpicture}
\hspace{1em}
\begin{tikzpicture}[scale=0.3]
\node[circle, draw, fill=black, inner sep=1.3pt] (A1) at (3,0) {};
\node[circle, draw, fill=black, inner sep=1.3pt] (A2) at (3,2) {};
\node[circle, draw, fill=black, inner sep=1.3pt] (A3) at (3,4) {};
\node[circle, draw, fill=black, inner sep=1.3pt] (B1) at (0,0) {};
\node[circle, draw, fill=black, inner sep=1.3pt] (B2) at (0,2) {};
\node[circle, draw, fill=black, inner sep=1.3pt] (B3) at (0,4) {};
\node (J) at (1.5,-1) {$J_{20}$};
\draw[blue, dashed, thick] (A1) -- (B1);
\draw[blue, dashed, thick] (A1) -- (B2);
\draw[blue, dashed, thick] (A1) -- (B3);
\draw[blue, dashed, thick] (A2) -- (B1);
\draw[blue, dashed, thick] (A2) -- (B2);
\draw[red, thick] (A2) -- (B3);
\draw[red, thick] (A3) -- (B1);
\draw[red, thick] (A3) -- (B2);
\draw[blue, dashed, thick] (A3) -- (B3);
\end{tikzpicture}
\hspace{1em}
\begin{tikzpicture}[scale=0.3]
\node[circle, draw, fill=black, inner sep=1.3pt] (A1) at (3,0) {};
\node[circle, draw, fill=black, inner sep=1.3pt] (A2) at (3,2) {};
\node[circle, draw, fill=black, inner sep=1.3pt] (A3) at (3,4) {};
\node[circle, draw, fill=black, inner sep=1.3pt] (B1) at (0,0) {};
\node[circle, draw, fill=black, inner sep=1.3pt] (B2) at (0,2) {};
\node[circle, draw, fill=black, inner sep=1.3pt] (B3) at (0,4) {};
\node (J) at (1.5,-1) {$J_{21}$};
\draw[blue, dashed, thick] (A1) -- (B1);
\draw[blue, dashed, thick] (A1) -- (B2);
\draw[blue, dashed, thick] (A1) -- (B3);
\draw[blue, dashed, thick] (A2) -- (B1);
\draw[blue, dashed, thick] (A2) -- (B2);
\draw[red, thick] (A2) -- (B3);
\draw[blue, dashed, thick] (A3) -- (B1);
\draw[red, thick] (A3) -- (B2);
\draw[red, thick] (A3) -- (B3);
\end{tikzpicture}
\hspace{1em}
\begin{tikzpicture}[scale=0.3]
\node[circle, draw, fill=black, inner sep=1.3pt] (A1) at (3,0) {};
\node[circle, draw, fill=black, inner sep=1.3pt] (A2) at (3,2) {};
\node[circle, draw, fill=black, inner sep=1.3pt] (A3) at (3,4) {};
\node[circle, draw, fill=black, inner sep=1.3pt] (B1) at (0,0) {};
\node[circle, draw, fill=black, inner sep=1.3pt] (B2) at (0,2) {};
\node[circle, draw, fill=black, inner sep=1.3pt] (B3) at (0,4) {};
\node (J) at (1.5,-1) {$J_{22}$};
\draw[blue, dashed, thick] (A1) -- (B1);
\draw[blue, dashed, thick] (A1) -- (B2);
\draw[blue, dashed, thick] (A1) -- (B3);
\draw[blue, dashed, thick] (A2) -- (B1);
\draw[blue, dashed, thick] (A2) -- (B2);
\draw[red, thick] (A2) -- (B3);
\draw[blue, dashed, thick] (A3) -- (B1);
\draw[red, thick] (A3) -- (B2);
\draw[blue, dashed, thick] (A3) -- (B3);
\end{tikzpicture}
\hspace{1em}
\begin{tikzpicture}[scale=0.3]
\node[circle, draw, fill=black, inner sep=1.3pt] (A1) at (3,0) {};
\node[circle, draw, fill=black, inner sep=1.3pt] (A2) at (3,2) {};
\node[circle, draw, fill=black, inner sep=1.3pt] (A3) at (3,4) {};
\node[circle, draw, fill=black, inner sep=1.3pt] (B1) at (0,0) {};
\node[circle, draw, fill=black, inner sep=1.3pt] (B2) at (0,2) {};
\node[circle, draw, fill=black, inner sep=1.3pt] (B3) at (0,4) {};
\node (J) at (1.5,-1) {$J_{23}$};
\draw[blue, dashed, thick] (A1) -- (B1);
\draw[blue, dashed, thick] (A1) -- (B2);
\draw[blue, dashed, thick] (A1) -- (B3);
\draw[blue, dashed, thick] (A2) -- (B1);
\draw[blue, dashed, thick] (A2) -- (B2);
\draw[blue, dashed, thick] (A2) -- (B3);
\draw[red, thick] (A3) -- (B1);
\draw[red, thick] (A3) -- (B2);
\draw[red, thick] (A3) -- (B3);
\end{tikzpicture}
\hspace{1em}
\begin{tikzpicture}[scale=0.3]
\node[circle, draw, fill=black, inner sep=1.3pt] (A1) at (3,0) {};
\node[circle, draw, fill=black, inner sep=1.3pt] (A2) at (3,2) {};
\node[circle, draw, fill=black, inner sep=1.3pt] (A3) at (3,4) {};
\node[circle, draw, fill=black, inner sep=1.3pt] (B1) at (0,0) {};
\node[circle, draw, fill=black, inner sep=1.3pt] (B2) at (0,2) {};
\node[circle, draw, fill=black, inner sep=1.3pt] (B3) at (0,4) {};
\node (J) at (1.5,-1) {$J_{24}$};
\draw[blue, dashed, thick] (A1) -- (B1);
\draw[blue, dashed, thick] (A1) -- (B2);
\draw[blue, dashed, thick] (A1) -- (B3);
\draw[blue, dashed, thick] (A2) -- (B1);
\draw[blue, dashed, thick] (A2) -- (B2);
\draw[blue, dashed, thick] (A2) -- (B3);
\draw[blue, dashed, thick] (A3) -- (B1);
\draw[red, thick] (A3) -- (B2);
\draw[red, thick] (A3) -- (B3);
\end{tikzpicture}
\hspace{1em}
\begin{tikzpicture}[scale=0.3]
\node[circle, draw, fill=black, inner sep=1.3pt] (A1) at (3,0) {};
\node[circle, draw, fill=black, inner sep=1.3pt] (A2) at (3,2) {};
\node[circle, draw, fill=black, inner sep=1.3pt] (A3) at (3,4) {};
\node[circle, draw, fill=black, inner sep=1.3pt] (B1) at (0,0) {};
\node[circle, draw, fill=black, inner sep=1.3pt] (B2) at (0,2) {};
\node[circle, draw, fill=black, inner sep=1.3pt] (B3) at (0,4) {};
\node (J) at (1.5,-1) {$J_{25}$};
\draw[blue, dashed, thick] (A1) -- (B1);
\draw[blue, dashed, thick] (A1) -- (B2);
\draw[blue, dashed, thick] (A1) -- (B3);
\draw[blue, dashed, thick] (A2) -- (B1);
\draw[blue, dashed, thick] (A2) -- (B2);
\draw[blue, dashed, thick] (A2) -- (B3);
\draw[blue, dashed, thick] (A3) -- (B1);
\draw[blue, dashed, thick] (A3) -- (B2);
\draw[red, thick] (A3) -- (B3);
\end{tikzpicture}
\hspace{1em}
\begin{tikzpicture}[scale=0.3]
\node[circle, draw, fill=black, inner sep=1.3pt] (A1) at (3,0) {};
\node[circle, draw, fill=black, inner sep=1.3pt] (A2) at (3,2) {};
\node[circle, draw, fill=black, inner sep=1.3pt] (A3) at (3,4) {};
\node[circle, draw, fill=black, inner sep=1.3pt] (B1) at (0,0) {};
\node[circle, draw, fill=black, inner sep=1.3pt] (B2) at (0,2) {};
\node[circle, draw, fill=black, inner sep=1.3pt] (B3) at (0,4) {};
\node (J) at (1.5,-1) {$J_{26}$};
\draw[blue, dashed, thick] (A1) -- (B1);
\draw[blue, dashed, thick] (A1) -- (B2);
\draw[blue, dashed, thick] (A1) -- (B3);
\draw[blue, dashed, thick] (A2) -- (B1);
\draw[blue, dashed, thick] (A2) -- (B2);
\draw[blue, dashed, thick] (A2) -- (B3);
\draw[blue, dashed, thick] (A3) -- (B1);
\draw[blue, dashed, thick] (A3) -- (B2);
\draw[blue, dashed, thick] (A3) -- (B3);
\end{tikzpicture}

    \caption{The 26 edge-colourings of $K_{3,3}$ up to isomorphism. Red edges are drawn as solid lines and blue edges are drawn as dashed lines.}
    \label{fig:K33}
\end{figure}
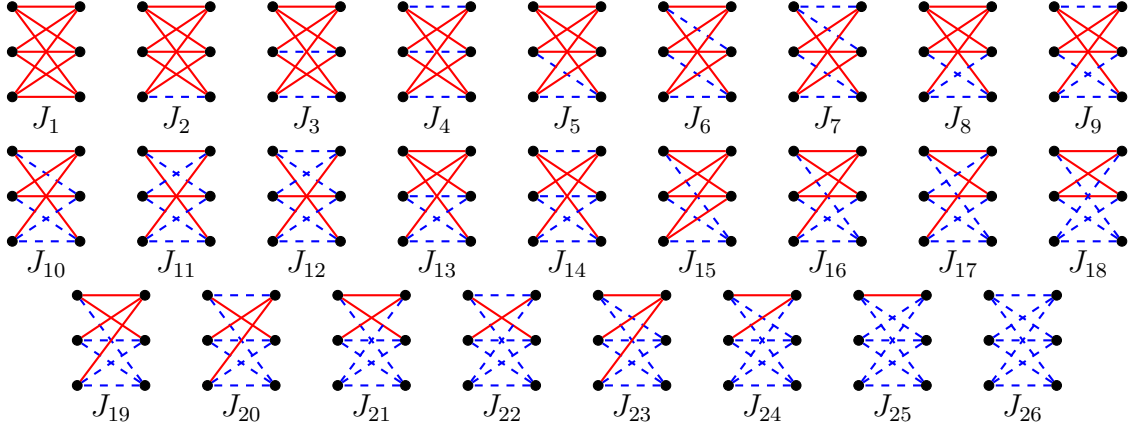

Of course, among the $512$ distinct edge-colourings of $K_{3,3}$, many of them are isomorphic to one another. Up to isomorphism, there are 26 different red/blue edge colourings of the complete bipartite graph $K_{3,3}$. Let $J_1,\dots,J_{26}$ be the corresponding 26 edge-coloured graphs, as depicted in Figure~\ref{fig:K33}, where $V(J_\ell)=V(K_{3,3})$ for all $1\leq\ell\leq26$. 

For an edge-coloured graph $H$, an \emph{automorphism} of $H$ is a bijective homomorphism from $H$ to itself. Let $\aut(H)$ be the number of automorphisms of $H$. For each $J_\ell$, the number of indices $j\in\{1,\dots,512\}$ such that $J_\ell$ is isomorphic to $H_j$ is precisely $\frac{2\cdot 3!\cdot 3!}{\aut(J_\ell)} = \frac{72}{\aut(J_\ell)}$. Therefore, if we let $d(J_\ell,G):=\frac{72}{\aut(J_\ell)}t_{\inj}(J_\ell,G)$ for each $1\leq \ell\leq 26$, then \eqref{eq:sumToOne} can be rewritten as
\begin{equation}
\label{eq:sumToOne2}
    \sum_{\ell=1}^{26}d(J_\ell,G) = 1.
\end{equation}
A \emph{homomorphism} between two graphs without edge colours is simply a map from the vertex set of the first graph to the vertex set of the second with the property that adjacent pairs of vertices are mapped to adjacent pairs of vertices. We view the graph $C_6$ and the edge-coloured graph $C_6^A$ as having the same vertex set. For each $1\leq \ell\leq 26$, let $t_{\bip}(C_6^A,J_\ell)$ be the probability that a uniformly random injective homomorphism $f$ from $C_6$ to $K_{3,3}$ is a homomorphism from $C_6^A$ to $J_\ell$; we define $t_{\bip}(C_6^A,H_i)$ analogously for $1\leq i\leq 512$. The next proposition relates $t_{\inj}(C_6^A,G)$ to the quantities $t_{\bip}(C_6^A,J_\ell)$ and $d(J_\ell,G)$ for $1\leq \ell\leq 26$. 

\begin{proposition}
\label{prop:doubleCount}
For any edge-coloured clique $G$,
\[t_{\inj}(C_6^A,G)=\sum_{\ell=1}^{26}t_{\bip}(C_6^A,J_\ell)\cdot d(J_\ell,G).\]
\end{proposition}

\begin{proof}
If $v(G)<6$, then both sides of the equation in the proposition are zero and so we may assume that $v(G)\geq6$. Let $f:V(K_{3,3})\to V(G)$ be a uniformly random injective function and let $g:V(C_6)\to V(K_{3,3})$ be a uniformly random injective homomorphism from $C_6$ to $K_{3,3}$. Then $f\circ g$ is nothing more than a uniformly random injective function from $V(C_6^A)$ to $V(G)$. Let $E$ be the event that $f\circ g$ is a homomorphism from $C_6^A$ to $G$ and, for each $1\leq i\leq 512$, let $F_i$ be the event that $f$ is a homomorphism from $H_i$ to $G$. By definition, we have $t_{\inj}(C_6^A,G)=\mathbb{P}(E)$. Also, since $G$ is an edge-coloured clique, the events $F_1,\dots,F_{512}$ are disjoint and satisfy $\sum_{i=1}^{512}\mathbb{P}(F_i)=1$. Therefore, 
\[t_{\inj}(C_6^A,G)=\sum_{i=1}^{512}\mathbb{P}(E\cap F_i)=\sum_{i=1}^{512}\mathbb{P}(E\mid F_i)\mathbb{P}(F_i) = \sum_{i=1}^{512}t_{\bip}(C_6^A,H_i)\cdot t_{\inj}(H_i,G).\]
The result follows by collecting terms based on isomorphism classes. 
\end{proof}

The quantities $t_{\bip}(C_6^A,J_\ell)$ for $1\leq \ell\leq 26$ can all be easily computed. As it turns out, 
\[t_{\bip}(C_6^A,J_4)=\frac{1}{6},\qquad t_{\bip}(C_6^A,J_9)=\frac{1}{12},\qquad t_{\bip}(C_6^A,J_{11})=\frac{1}{12},\qquad t_{\bip}(C_6^A,J_{12})=\frac{1}{6}\]
and $t_{\bip}(C_6^A,J_\ell)=0$ for all other indices $\ell$. Thus, by Proposition~\ref{prop:doubleCount},
\begin{equation}
\label{eq:doubleCount}
t_{\inj}(C_6^A,G)=\frac{1}{6}d(J_4,G) + \frac{1}{12}d(J_9,G) +\frac{1}{12}d(J_{11},G) +\frac{1}{6}d(J_{12},G).
\end{equation}

\begin{figure}[htbp]
    \centering
\begin{tikzpicture}[scale=0.3]
\node[rectangle, draw, fill=black, inner sep=1.3pt] (A1) at (3,0) {};
\node[circle, draw, fill=black, inner sep=1.3pt] (A2) at (3,2) {};
\node[rectangle, draw, fill=black, inner sep=1.3pt] (B1) at (0,0) {};
\node[circle, draw, fill=black, inner sep=1.3pt] (B2) at (0,2) {};
\node (R) at (1.5,-1) {$R_{1}$};
\draw[red, thick] (A1) -- (B1);
\draw[red, thick] (A1) -- (B2);
\draw[red, thick] (A2) -- (B1);
\draw[red, thick] (A2) -- (B2);

\node (one) at (-0.7,0) {\footnotesize 1};
\node (two) at (3.7,0) {\footnotesize 2};
\end{tikzpicture}
\hspace{0.1em}
\begin{tikzpicture}[scale=0.3]
\node[rectangle, draw, fill=black, inner sep=1.3pt] (A1) at (3,0) {};
\node[circle, draw, fill=black, inner sep=1.3pt] (A2) at (3,2) {};
\node[rectangle, draw, fill=black, inner sep=1.3pt] (B1) at (0,0) {};
\node[circle, draw, fill=black, inner sep=1.3pt] (B2) at (0,2) {};
\node (R) at (1.5,-1) {$R_{2}$};
\draw[red, thick] (A1) -- (B1);
\draw[red, thick] (A1) -- (B2);
\draw[red, thick] (A2) -- (B1);
\draw[blue, dashed, thick] (A2) -- (B2);

\node (one) at (-0.7,0) {\footnotesize 1};
\node (two) at (3.7,0) {\footnotesize 2};
\end{tikzpicture}
\hspace{0.1em}
\begin{tikzpicture}[scale=0.3]
\node[rectangle, draw, fill=black, inner sep=1.3pt] (A1) at (3,0) {};
\node[circle, draw, fill=black, inner sep=1.3pt] (A2) at (3,2) {};
\node[rectangle, draw, fill=black, inner sep=1.3pt] (B1) at (0,0) {};
\node[circle, draw, fill=black, inner sep=1.3pt] (B2) at (0,2) {};
\node (R) at (1.5,-1) {$R_{3}$};
\draw[red, thick] (A1) -- (B1);
\draw[red, thick] (A1) -- (B2);
\draw[blue, dashed, thick] (A2) -- (B1);
\draw[red, thick] (A2) -- (B2);

\node (one) at (-0.7,0) {\footnotesize 1};
\node (two) at (3.7,0) {\footnotesize 2};
\end{tikzpicture}
\hspace{0.1em}
\begin{tikzpicture}[scale=0.3]
\node[rectangle, draw, fill=black, inner sep=1.3pt] (A1) at (3,0) {};
\node[circle, draw, fill=black, inner sep=1.3pt] (A2) at (3,2) {};
\node[rectangle, draw, fill=black, inner sep=1.3pt] (B1) at (0,0) {};
\node[circle, draw, fill=black, inner sep=1.3pt] (B2) at (0,2) {};
\node (R) at (1.5,-1) {$R_{4}$};
\draw[red, thick] (A1) -- (B1);
\draw[red, thick] (A1) -- (B2);
\draw[blue, dashed, thick] (A2) -- (B1);
\draw[blue, dashed, thick] (A2) -- (B2);

\node (one) at (-0.7,0) {\footnotesize 1};
\node (two) at (3.7,0) {\footnotesize 2};
\end{tikzpicture}
\hspace{0.1em}
\begin{tikzpicture}[scale=0.3]
\node[rectangle, draw, fill=black, inner sep=1.3pt] (A1) at (3,0) {};
\node[circle, draw, fill=black, inner sep=1.3pt] (A2) at (3,2) {};
\node[rectangle, draw, fill=black, inner sep=1.3pt] (B1) at (0,0) {};
\node[circle, draw, fill=black, inner sep=1.3pt] (B2) at (0,2) {};
\node (R) at (1.5,-1) {$R_{5}$};
\draw[red, thick] (A1) -- (B1);
\draw[blue, dashed, thick] (A1) -- (B2);
\draw[red, thick] (A2) -- (B1);
\draw[red, thick] (A2) -- (B2);

\node (one) at (-0.7,0) {\footnotesize 1};
\node (two) at (3.7,0) {\footnotesize 2};
\end{tikzpicture}
\hspace{0.1em}
\begin{tikzpicture}[scale=0.3]
\node[rectangle, draw, fill=black, inner sep=1.3pt] (A1) at (3,0) {};
\node[circle, draw, fill=black, inner sep=1.3pt] (A2) at (3,2) {};
\node[rectangle, draw, fill=black, inner sep=1.3pt] (B1) at (0,0) {};
\node[circle, draw, fill=black, inner sep=1.3pt] (B2) at (0,2) {};
\node (R) at (1.5,-1) {$R_{6}$};
\draw[red, thick] (A1) -- (B1);
\draw[blue, dashed, thick] (A1) -- (B2);
\draw[red, thick] (A2) -- (B1);
\draw[blue, dashed, thick] (A2) -- (B2);

\node (one) at (-0.7,0) {\footnotesize 1};
\node (two) at (3.7,0) {\footnotesize 2};
\end{tikzpicture}
\hspace{0.1em}
\begin{tikzpicture}[scale=0.3]
\node[rectangle, draw, fill=black, inner sep=1.3pt] (A1) at (3,0) {};
\node[circle, draw, fill=black, inner sep=1.3pt] (A2) at (3,2) {};
\node[rectangle, draw, fill=black, inner sep=1.3pt] (B1) at (0,0) {};
\node[circle, draw, fill=black, inner sep=1.3pt] (B2) at (0,2) {};
\node (R) at (1.5,-1) {$R_{7}$};
\draw[red, thick] (A1) -- (B1);
\draw[blue, dashed, thick] (A1) -- (B2);
\draw[blue, dashed, thick] (A2) -- (B1);
\draw[red, thick] (A2) -- (B2);

\node (one) at (-0.7,0) {\footnotesize 1};
\node (two) at (3.7,0) {\footnotesize 2};
\end{tikzpicture}
\hspace{0.1em}
\begin{tikzpicture}[scale=0.3]
\node[rectangle, draw, fill=black, inner sep=1.3pt] (A1) at (3,0) {};
\node[circle, draw, fill=black, inner sep=1.3pt] (A2) at (3,2) {};
\node[rectangle, draw, fill=black, inner sep=1.3pt] (B1) at (0,0) {};
\node[circle, draw, fill=black, inner sep=1.3pt] (B2) at (0,2) {};
\node (R) at (1.5,-1) {$R_{8}$};
\draw[red, thick] (A1) -- (B1);
\draw[blue, dashed, thick] (A1) -- (B2);
\draw[blue, dashed, thick] (A2) -- (B1);
\draw[blue, dashed, thick] (A2) -- (B2);

\node (one) at (-0.7,0) {\footnotesize 1};
\node (two) at (3.7,0) {\footnotesize 2};
\end{tikzpicture}

\begin{tikzpicture}[scale=0.3]
\node[rectangle, draw, fill=black, inner sep=1.3pt] (A1) at (3,0) {};
\node[circle, draw, fill=black, inner sep=1.3pt] (A2) at (3,2) {};
\node[rectangle, draw, fill=black, inner sep=1.3pt] (B1) at (0,0) {};
\node[circle, draw, fill=black, inner sep=1.3pt] (B2) at (0,2) {};
\node (R) at (1.5,-1) {$B_{1}$};
\draw[blue, dashed, thick] (A1) -- (B1);
\draw[blue, dashed, thick] (A1) -- (B2);
\draw[blue, dashed, thick] (A2) -- (B1);
\draw[blue, dashed, thick] (A2) -- (B2);

\node (one) at (-0.7,0) {\footnotesize 1};
\node (two) at (3.7,0) {\footnotesize 2};
\end{tikzpicture}
\hspace{0.1em}
\begin{tikzpicture}[scale=0.3]
\node[rectangle, draw, fill=black, inner sep=1.3pt] (A1) at (3,0) {};
\node[circle, draw, fill=black, inner sep=1.3pt] (A2) at (3,2) {};
\node[rectangle, draw, fill=black, inner sep=1.3pt] (B1) at (0,0) {};
\node[circle, draw, fill=black, inner sep=1.3pt] (B2) at (0,2) {};
\node (R) at (1.5,-1) {$B_{2}$};
\draw[blue, dashed, thick] (A1) -- (B1);
\draw[blue, dashed, thick] (A1) -- (B2);
\draw[blue, dashed, thick] (A2) -- (B1);
\draw[red, thick] (A2) -- (B2);

\node (one) at (-0.7,0) {\footnotesize 1};
\node (two) at (3.7,0) {\footnotesize 2};
\end{tikzpicture}
\hspace{0.1em}
\begin{tikzpicture}[scale=0.3]
\node[rectangle, draw, fill=black, inner sep=1.3pt] (A1) at (3,0) {};
\node[circle, draw, fill=black, inner sep=1.3pt] (A2) at (3,2) {};
\node[rectangle, draw, fill=black, inner sep=1.3pt] (B1) at (0,0) {};
\node[circle, draw, fill=black, inner sep=1.3pt] (B2) at (0,2) {};
\node (R) at (1.5,-1) {$B_{3}$};
\draw[blue, dashed, thick] (A1) -- (B1);
\draw[blue, dashed, thick] (A1) -- (B2);
\draw[red, thick] (A2) -- (B1);
\draw[blue, dashed, thick] (A2) -- (B2);

\node (one) at (-0.7,0) {\footnotesize 1};
\node (two) at (3.7,0) {\footnotesize 2};
\end{tikzpicture}
\hspace{0.1em}
\begin{tikzpicture}[scale=0.3]
\node[rectangle, draw, fill=black, inner sep=1.3pt] (A1) at (3,0) {};
\node[circle, draw, fill=black, inner sep=1.3pt] (A2) at (3,2) {};
\node[rectangle, draw, fill=black, inner sep=1.3pt] (B1) at (0,0) {};
\node[circle, draw, fill=black, inner sep=1.3pt] (B2) at (0,2) {};
\node (R) at (1.5,-1) {$B_{4}$};
\draw[blue, dashed, thick] (A1) -- (B1);
\draw[blue, dashed, thick] (A1) -- (B2);
\draw[red, thick] (A2) -- (B1);
\draw[red, thick] (A2) -- (B2);

\node (one) at (-0.7,0) {\footnotesize 1};
\node (two) at (3.7,0) {\footnotesize 2};
\end{tikzpicture}
\hspace{0.1em}
\begin{tikzpicture}[scale=0.3]
\node[rectangle, draw, fill=black, inner sep=1.3pt] (A1) at (3,0) {};
\node[circle, draw, fill=black, inner sep=1.3pt] (A2) at (3,2) {};
\node[rectangle, draw, fill=black, inner sep=1.3pt] (B1) at (0,0) {};
\node[circle, draw, fill=black, inner sep=1.3pt] (B2) at (0,2) {};
\node (R) at (1.5,-1) {$B_{5}$};
\draw[blue, dashed, thick] (A1) -- (B1);
\draw[red, thick] (A1) -- (B2);
\draw[blue, dashed, thick] (A2) -- (B1);
\draw[blue, dashed, thick] (A2) -- (B2);

\node (one) at (-0.7,0) {\footnotesize 1};
\node (two) at (3.7,0) {\footnotesize 2};
\end{tikzpicture}
\hspace{0.1em}
\begin{tikzpicture}[scale=0.3]
\node[rectangle, draw, fill=black, inner sep=1.3pt] (A1) at (3,0) {};
\node[circle, draw, fill=black, inner sep=1.3pt] (A2) at (3,2) {};
\node[rectangle, draw, fill=black, inner sep=1.3pt] (B1) at (0,0) {};
\node[circle, draw, fill=black, inner sep=1.3pt] (B2) at (0,2) {};
\node (R) at (1.5,-1) {$B_{6}$};
\draw[blue, dashed, thick] (A1) -- (B1);
\draw[red, thick] (A1) -- (B2);
\draw[blue, dashed, thick] (A2) -- (B1);
\draw[red, thick] (A2) -- (B2);

\node (one) at (-0.7,0) {\footnotesize 1};
\node (two) at (3.7,0) {\footnotesize 2};
\end{tikzpicture}
\hspace{0.1em}
\begin{tikzpicture}[scale=0.3]
\node[rectangle, draw, fill=black, inner sep=1.3pt] (A1) at (3,0) {};
\node[circle, draw, fill=black, inner sep=1.3pt] (A2) at (3,2) {};
\node[rectangle, draw, fill=black, inner sep=1.3pt] (B1) at (0,0) {};
\node[circle, draw, fill=black, inner sep=1.3pt] (B2) at (0,2) {};
\node (R) at (1.5,-1) {$B_{7}$};
\draw[blue, dashed, thick] (A1) -- (B1);
\draw[red, thick] (A1) -- (B2);
\draw[red, thick] (A2) -- (B1);
\draw[blue, dashed, thick] (A2) -- (B2);

\node (one) at (-0.7,0) {\footnotesize 1};
\node (two) at (3.7,0) {\footnotesize 2};
\end{tikzpicture}
\hspace{0.1em}
\begin{tikzpicture}[scale=0.3]
\node[rectangle, draw, fill=black, inner sep=1.3pt] (A1) at (3,0) {};
\node[circle, draw, fill=black, inner sep=1.3pt] (A2) at (3,2) {};
\node[rectangle, draw, fill=black, inner sep=1.3pt] (B1) at (0,0) {};
\node[circle, draw, fill=black, inner sep=1.3pt] (B2) at (0,2) {};
\node (R) at (1.5,-1) {$B_{8}$};
\draw[blue, dashed, thick] (A1) -- (B1);
\draw[red, thick] (A1) -- (B2);
\draw[red, thick] (A2) -- (B1);
\draw[red, thick] (A2) -- (B2);

\node (one) at (-0.7,0) {\footnotesize 1};
\node (two) at (3.7,0) {\footnotesize 2};
\end{tikzpicture}
    \caption{The flags $R_1,\dots,R_8$ and $B_1,\dots,B_8$. The roots are depicted with square nodes. }
    \label{fig:flags}
\end{figure}
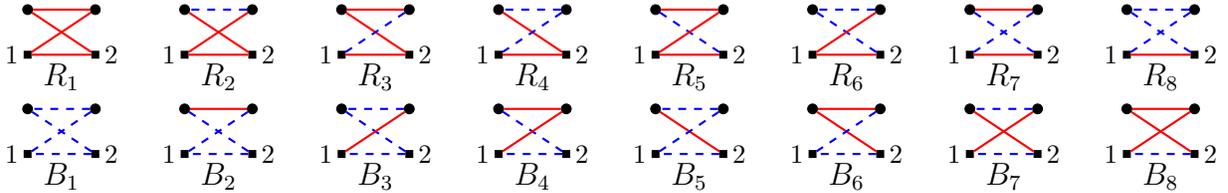

Now, consider the $16$ graphs $R_1,\dots,R_8$ and $B_1,\dots,B_8$ in Figure~\ref{fig:flags}. We refer to these graphs as \emph{flags} and the vertices labelled $1$ and $2$ in each flag are referred to as the \emph{roots} of the flag. Given a pair $u,v$ of distinct vertices of $G$ and a graph $F$ with two root vertices labelled $1$ and $2$, let $\hom_{\inj}(F,G; u,v)$ be the number of injective homomorphisms $f:V(F)\to V(G)$ with $f(1)=u$ and $f(2)=v$. 

\begin{proof}[Proof of Theorem~\ref{th:injective}]
For each pair $u,v$ of distinct vertices of $G$, let $\vec{x_{u,v}}$ and $\vec{y_{u,v}}$ be the vectors $\vec{x_{u,v}}:=[\hom_{\inj}(R_i,G;u,v): 1\leq i\leq 8]^T$ and $\vec{y_{u,v}}:=[\hom_{\inj}(B_i,G;u,v): 1\leq i\leq 8]^T$. We define an $8\times 8$ matrix $A$ as follows:
\[A := \frac{1}{128}\begin{bmatrix}
2& -6& -2& -3& 1& -3& 5& 6\\
-6& 58& -3& 12& -6& 12& -47& -20\\
-2& -3& 56& -14& -47& 11& 4& -5\\
-3& 12& -14& 12 &10& 2 &-10& -9\\
1& -6& -47& 10& 56& -14& 2 &-2\\
-3& 12& 11& 2 &-14& 12& -10& -10\\
5& -47& 4& -10& 2& -10& 44& 12\\
6& -20& -5& -9& -2& -10& 12& 28\end{bmatrix}.\]
The matrix $A$ is positive semi-definite (i.e. all of its eigenvalues are non-negative). Therefore, $\vec{x_{u,v}}^TA\vec{x_{u,v}}\geq0$ and $\vec{y_{u,v}}^TA\vec{y_{u,v}}\geq0$ for every pair $u,v$ of distinct vertices of $G$. So, by \eqref{eq:doubleCount}, $t_{\inj}(C_6^A,G)$ is at most
\begin{equation}\label{eq:flagged}\scalemath{0.85}{\frac{1}{6}d(J_4,G) + \frac{1}{12}d(J_9,G) +\frac{1}{12}d(J_{11},G) +\frac{1}{6}d(J_{12},G)+\frac{(n-6)!}{n!}\left(\sum_{u,v}\vec{x_{u,v}}^TA\vec{x_{u,v}} + \sum_{u,v}\vec{y_{u,v}}^TA\vec{y_{u,v}}\right).}\end{equation}
By \eqref{eq:sumToOne2}, it suffices to show that \eqref{eq:flagged} evaluates to $(1/2)^6\sum_{\ell=1}^{26}d(J_\ell,G) + o(1)$. In order to do this, we need to do some ``double counting'' to express the summations in \eqref{eq:flagged} in terms of the quantities $d(J_\ell,G)$ and a small error term. We start by observing that
\[\frac{(n-6)!}{n!}\sum_{u,v}\vec{x_{u,v}}^TA\vec{x_{u,v}} =\sum_{i=1}^8\sum_{j=1}^8A(i,j)\frac{(n-6)!}{n!}\sum_{u,v}\hom_{\inj}(R_i,G;u,v)\hom_{\inj}(R_j,G;u,v)\]
where $A(i,j)$ is the entry of $A$ on the $i$th row and $j$th column and, similarly,
\[\frac{(n-6)!}{n!}\sum_{u,v}\vec{y_{u,v}}^TA\vec{y_{u,v}} =\sum_{i=1}^8\sum_{j=1}^8A(i,j)\frac{(n-6)!}{n!}\sum_{u,v}\hom_{\inj}(B_i,G;u,v)\hom_{\inj}(B_j,G;u,v).\]
For distinct vertices $u$ and $v$, the quantity $\hom_{\inj}(R_i,G;u,v)\hom_{\inj}(R_j,G;u,v)$ counts pairs $(f_i,f_j)$ where $f_i$ is an injective homomorphism from $R_i$ to $G$ and $f_j$ is an injective homomorphism from $R_j$ to $G$ such that $f_i(1)=f_j(1)=u$ and $f_i(2)=f_j(2)=v$. Let $R_i\cdot R_j$ denote the $6$-vertex graph obtained from the disjoint union of $R_i$ and $R_j$ by identifying each root of $R_i$ with the corresponding root of $R_j$. Then $\hom_{\inj}(R_i,G;u,v)\hom_{\inj}(R_j,G;u,v)$ is equal to $\hom_{\inj}(R_i\cdot R_j,G; u,v) + O(n^3)$ where the $O(n^3)$ term accounts for the pairs $(f_i,f_j)$ counted by $\hom_{\inj}(R_i,G;u,v)\hom_{\inj}(R_j,G;u,v)$ such that $f_i(w)=f_j(z)$ for some $w\in V(R_i)\setminus\{1,2\}$ and $z\in V(R_j)\setminus\{1,2\}$. Similarly, $\hom_{\inj}(B_i,G;u,v)\hom_{\inj}(B_j,G;u,v)$ can be written as $\hom_{\inj}(B_i\cdot B_j,G; u,v) + O(n^3)$ where the graph $B_i\cdot B_j$ is defined analogously to $R_i\cdot R_j$. The contribution of the $O(n^3)$ terms to \eqref{eq:flagged} is $O(1/n) = o(1)$ and so, up to a $o(1)$ term, we can replace $\hom_{\inj}(R_i,G;u,v)\hom_{\inj}(R_j,G;u,v)$ with $\hom_{\inj}(R_i\cdot R_j,G; u,v)$ and $\hom_{\inj}(B_i,G;u,v)\hom_{\inj}(B_j,G;u,v)$ with $\hom_{\inj}(B_i\cdot B_j,G; u,v)$. 

Observe that $\frac{(n-6)!}{n!}\sum_{u,v}\hom_{\inj}(R_i\cdot R_j,G; u,v) = t_{\inj}(R_i\cdot R_j,G)$. Combining with~\eqref{eq:flagged}, we obtain that $t_{\inj}(C_6^A,G)$ is at most
\begin{equation}\label{C6upperbound}
\begin{aligned}
&\frac{1}{6}d(J_4,G) + \frac{1}{12}d(J_9,G) +\frac{1}{12}d(J_{11},G) +\frac{1}{6}d(J_{12},G)\\
&+\sum_{i=1}^8\sum_{j=1}^8A(i,j)\left(t_{\inj}(R_i\cdot R_j,G) +t_{\inj}(B_i\cdot B_j,G) \right)+o(1).
\end{aligned}
\end{equation}

The following claim will complete the proof of Theorem~\ref{th:injective}. Since the proof of this claim is entirely numerical and a bit long, we decide to put its proof in Appendix A.

\begin{claim}\label{numbercheck}
It holds that
\begin{equation}
\begin{aligned}
&\frac{1}{6}d(J_4,G) + \frac{1}{12}d(J_9,G) +\frac{1}{12}d(J_{11},G) +\frac{1}{6}d(J_{12},G)\\
&+\sum_{i=1}^8\sum_{j=1}^8A(i,j)\left(t_{\inj}(R_i\cdot R_j,G) +t_{\inj}(B_i\cdot B_j,G) \right)=\left(\frac{1}{2}\right)^6.
\end{aligned}
\end{equation}
\end{claim}

By this claim, we have that~\eqref{C6upperbound} is equal to $(1/2)^6+o(1)$, which implies $t_{\inj}(C_6^A,G)\leq (1/2)^6 + o(1)$. 
\end{proof}

Say that a sequence $G_1,G_2,\dots$ of edge-coloured cliques is \emph{quasirandom} if, for every edge-coloured clique $H$,
\[\lim_{n\to\infty}t(H,G_n)=(1/2)^{e(H)}.\]
By analyzing the proof of Theorem~\ref{th:injective} further, one can get that a sequence $G_1,G_2,\dots$ of edge-coloured cliques satisfies $\lim_{n\to\infty}t(C_6^A,G_n)=(1/2)^6$ if and only if $G_1,G_2,\dots$ is quasirandom. Here is a brief summary of the argument. We observe that the kernel of the matrix $A$ is spanned by the all-ones vector of length 8. This implies that, for large $n$, nearly every pair of distinct vertices is connected by $\left(\frac{1}{8}+o(1)\right)n^2$ red paths of length three. As a result, we get that the homomorphism density of $C_6$ in the graph consisting of the red edges of $G$ is $(1/2)^6 + o(1)$. By standard results on quasirandomness of graphs, this implies that the graph consisting of the red edges is quasirandom, which also implies that the sequence $G_1,G_2,\dots$ is quasirandom. 

A limitation of the flag algebra method is that extremal problems involving graphs on more than a small number of vertices (perhaps somewhere between seven and eleven, depending on the problem and the skill level of the person writing the semi-definite programs) tend to be out of reach due to combinatorial explosion. We find it interesting that the flag algebra proof in this paper only needed to consider the 26 non-isomorphic colourings of $K_{3,3}$ as opposed to the 156 non-isomorphic colourings of $K_6$. Likewise, the flags that we used have missing edges which helped to reduce the number of flags that we required from 40 to 16. We wonder whether considering these types of ``simpler'' flag algebra proofs could, occasionally, allow one to extend the method to larger graphs than one can normally handle, or to obtain human-checkable proofs for theorems which currently require computer verification.

\subsection*{Acknowledgements}

This work was completed while Hao Chen was visiting the University of Victoria. Hao Chen would like to thank the China Scholarship Council for its support. Jonathan A. Noel was supported by NSERC Discovery Grant RGPIN-2021-02460.

\appendix
\section{Proof of Claim~\ref{numbercheck}}
In this section, we provide the proof of Claim~\ref{numbercheck}.

\begin{proof}[Proof of Claim~\ref{numbercheck}]
By following the same arguments as in Proposition~\ref{prop:doubleCount}, we obtain, for all $1\leq i,j\leq 8$,
\begin{equation*}
t_{\inj}(R_i\cdot R_j,G) = \sum_{\ell=1}^{26}t_{\bip}(R_i\cdot R_j,J_\ell)\cdot d(J_\ell,G)
\end{equation*}
and 
\begin{equation*}
t_{\inj}(B_i\cdot B_j,G) = \sum_{\ell=1}^{26}t_{\bip}(B_i\cdot B_j,J_\ell)\cdot d(J_\ell,G).
\end{equation*}

By computing the quantities $t_{\bip}(R_i\cdot R_j,J_\ell)$ and $t_{\bip}(B_i\cdot B_j,J_\ell)$ explicitly for all $1\leq i\leq j\leq 8$ and $1\leq \ell\leq 26$, we obtain the following:
\[t_{\inj}(R_{1}\cdot R_{1},G) = d(J_{1},G)+\frac{16}{72}d(J_{2},G)+\frac{4}{72}d(J_{3},G)\]
\[t_{\inj}(R_{1}\cdot R_{2},G) = \frac{8}{72}d(J_{2},G)+\frac{8}{72}d(J_{5},G)+\frac{2}{72}d(J_{8},G)\]
\[t_{\inj}(R_{1}\cdot R_{3},G) = \frac{8}{72}d(J_{2},G)+\frac{4}{72}d(J_{3},G)+\frac{4}{72}d(J_{5},G)+\frac{2}{72}d(J_{6},G)\]
\[t_{\inj}(R_{1}\cdot R_{4},G) = \frac{4}{72}d(J_{5},G)+\frac{2}{72}d(J_{8},G)+\frac{12}{72}d(J_{15},G)+\frac{2}{72}d(J_{16},G)\]
\[t_{\inj}(R_{1}\cdot R_{5},G) = \frac{8}{72}d(J_{2},G)+\frac{4}{72}d(J_{3},G)+\frac{4}{72}d(J_{5},G)+\frac{2}{72}d(J_{6},G)\]
\[t_{\inj}(R_{1}\cdot R_{6},G) = \frac{4}{72}d(J_{5},G)+\frac{2}{72}d(J_{8},G)+\frac{12}{72}d(J_{15},G)+\frac{2}{72}d(J_{16},G)\]
\[t_{\inj}(R_{1}\cdot R_{7},G) = \frac{4}{72}d(J_{3},G)+\frac{4}{72}d(J_{6},G)+\frac{8}{72}d(J_{7},G)\]
\[t_{\inj}(R_{1}\cdot R_{8},G) = \frac{2}{72}d(J_{8},G)+\frac{4}{72}d(J_{16},G)+\frac{8}{72}d(J_{18},G)\]
\[t_{\inj}(R_{2}\cdot R_{2},G) = \frac{4}{72}d(J_{3},G)+\frac{4}{72}d(J_{8},G)+\frac{8}{72}d(J_{13},G)\]
\[t_{\inj}(R_{2}\cdot R_{3},G) = \frac{4}{72}d(J_{3},G)+\frac{2}{72}d(J_{6},G)+\frac{2}{72}d(J_{8},G)+\frac{2}{72}d(J_{10},G)\]
\[t_{\inj}(R_{2}\cdot R_{4},G) = \frac{2}{72}d(J_{6},G)+\frac{2}{72}d(J_{10},G)+\frac{2}{72}d(J_{16},G)+\frac{2}{72}d(J_{19},G)\]
\[t_{\inj}(R_{2}\cdot R_{5},G) = \frac{4}{72}d(J_{3},G)+\frac{2}{72}d(J_{6},G)+\frac{2}{72}d(J_{8},G)+\frac{2}{72}d(J_{10},G)\]
\[t_{\inj}(R_{2}\cdot R_{6},G) = \frac{2}{72}d(J_{6},G)+\frac{2}{72}d(J_{10},G)+\frac{2}{72}d(J_{16},G)+\frac{2}{72}d(J_{19},G)\]
\[t_{\inj}(R_{2}\cdot R_{7},G) = \frac{12}{72}d(J_{4},G)+\frac{4}{72}d(J_{9},G)+\frac{2}{72}d(J_{11},G)\]
\[t_{\inj}(R_{2}\cdot R_{8},G) = \frac{2}{72}d(J_{9},G)+\frac{4}{72}d(J_{17},G)+\frac{2}{72}d(J_{21},G)\]
\[t_{\inj}(R_{3}\cdot R_{3},G) = \frac{4}{72}d(J_{5},G)+\frac{4}{72}d(J_{8},G)+\frac{2}{72}d(J_{10},G)\]
\[t_{\inj}(R_{3}\cdot R_{4},G) = \frac{2}{72}d(J_{8},G)+\frac{8}{72}d(J_{13},G)+\frac{2}{72}d(J_{16},G)+\frac{2}{72}d(J_{19},G)\]
\[t_{\inj}(R_{3}\cdot R_{5},G) = \frac{4}{72}d(J_{3},G)+\frac{12}{72}d(J_{4},G)+\frac{2}{72}d(J_{8},G)+\frac{2}{72}d(J_{9},G)\]
\[t_{\inj}(R_{3}\cdot R_{6},G) = \frac{2}{72}d(J_{6},G)+\frac{2}{72}d(J_{9},G)+\frac{2}{72}d(J_{16},G)+\frac{2}{72}d(J_{17},G)\]
\[t_{\inj}(R_{3}\cdot R_{7},G) = \frac{2}{72}d(J_{6},G)+\frac{2}{72}d(J_{9},G)+\frac{2}{72}d(J_{10},G)+\frac{2}{72}d(J_{11},G)\]
\[t_{\inj}(R_{3}\cdot R_{8},G) = \frac{2}{72}d(J_{10},G)+\frac{2}{72}d(J_{17},G)+\frac{2}{72}d(J_{19},G)+\frac{2}{72}d(J_{21},G)\]
\[t_{\inj}(R_{4}\cdot R_{4},G) = \frac{2}{72}d(J_{10},G)+\frac{4}{72}d(J_{19},G)+\frac{12}{72}d(J_{23},G)\]
\[t_{\inj}(R_{4}\cdot R_{5},G) = \frac{2}{72}d(J_{6},G)+\frac{2}{72}d(J_{9},G)+\frac{2}{72}d(J_{16},G)+\frac{2}{72}d(J_{17},G)\]
\[t_{\inj}(R_{4}\cdot R_{6},G) = \frac{8}{72}d(J_{7},G)+\frac{2}{72}d(J_{11},G)+\frac{8}{72}d(J_{18},G)+\frac{2}{72}d(J_{21},G)\]
\[t_{\inj}(R_{4}\cdot R_{7},G) = \frac{2}{72}d(J_{9},G)+\frac{8}{72}d(J_{14},G)+\frac{2}{72}d(J_{17},G)+\frac{2}{72}d(J_{20},G)\]
\[t_{\inj}(R_{4}\cdot R_{8},G) = \frac{2}{72}d(J_{11},G)+\frac{2}{72}d(J_{20},G)+\frac{2}{72}d(J_{21},G)+\frac{4}{72}d(J_{24},G)\]
\[t_{\inj}(R_{5}\cdot R_{5},G) = \frac{4}{72}d(J_{5},G)+\frac{4}{72}d(J_{8},G)+\frac{2}{72}d(J_{10},G)\]
\[t_{\inj}(R_{5}\cdot R_{6},G) = \frac{2}{72}d(J_{8},G)+\frac{8}{72}d(J_{13},G)+\frac{2}{72}d(J_{16},G)+\frac{2}{72}d(J_{19},G)\]
\[t_{\inj}(R_{5}\cdot R_{7},G) = \frac{2}{72}d(J_{6},G)+\frac{2}{72}d(J_{9},G)+\frac{2}{72}d(J_{10},G)+\frac{2}{72}d(J_{11},G)\]
\[t_{\inj}(R_{5}\cdot R_{8},G) = \frac{2}{72}d(J_{10},G)+\frac{2}{72}d(J_{17},G)+\frac{2}{72}d(J_{19},G)+\frac{2}{72}d(J_{21},G)\]
\[t_{\inj}(R_{6}\cdot R_{6},G) = \frac{2}{72}d(J_{10},G)+\frac{4}{72}d(J_{19},G)+\frac{12}{72}d(J_{23},G)\]
\[t_{\inj}(R_{6}\cdot R_{7},G) = \frac{2}{72}d(J_{9},G)+\frac{8}{72}d(J_{14},G)+\frac{2}{72}d(J_{17},G)+\frac{2}{72}d(J_{20},G)\]
\[t_{\inj}(R_{6}\cdot R_{8},G) = \frac{2}{72}d(J_{11},G)+\frac{2}{72}d(J_{20},G)+\frac{2}{72}d(J_{21},G)+\frac{4}{72}d(J_{24},G)\]
\[t_{\inj}(R_{7}\cdot R_{7},G) = \frac{8}{72}d(J_{7},G)+\frac{4}{72}d(J_{11},G)+\frac{12}{72}d(J_{12},G)\]
\[t_{\inj}(R_{7}\cdot R_{8},G) = \frac{2}{72}d(J_{11},G)+\frac{4}{72}d(J_{20},G)+\frac{4}{72}d(J_{22},G)\]
\[t_{\inj}(R_{8}\cdot R_{8},G) = \frac{12}{72}d(J_{12},G)+\frac{8}{72}d(J_{22},G)+\frac{8}{72}d(J_{25},G)\]
\[t_{\inj}(B_{1}\cdot B_{1},G) = \frac{4}{72}d(J_{22},G)+\frac{16}{72}d(J_{25},G)+d(J_{26},G)\]
\[t_{\inj}(B_{1}\cdot B_{2},G) = \frac{2}{72}d(J_{21},G)+\frac{8}{72}d(J_{24},G)+\frac{8}{72}d(J_{25},G)\]
\[t_{\inj}(B_{1}\cdot B_{3},G) = \frac{2}{72}d(J_{20},G)+\frac{4}{72}d(J_{22},G)+\frac{4}{72}d(J_{24},G)+\frac{8}{72}d(J_{25},G)\]
\[t_{\inj}(B_{1}\cdot B_{4},G) = \frac{2}{72}d(J_{19},G)+\frac{2}{72}d(J_{21},G)+\frac{12}{72}d(J_{23},G)+\frac{4}{72}d(J_{24},G)\]
\[t_{\inj}(B_{1}\cdot B_{5},G) = \frac{2}{72}d(J_{20},G)+\frac{4}{72}d(J_{22},G)+\frac{4}{72}d(J_{24},G)+\frac{8}{72}d(J_{25},G)\]
\[t_{\inj}(B_{1}\cdot B_{6},G) = \frac{2}{72}d(J_{19},G)+\frac{2}{72}d(J_{21},G)+\frac{12}{72}d(J_{23},G)+\frac{4}{72}d(J_{24},G)\]
\[t_{\inj}(B_{1}\cdot B_{7},G) = \frac{8}{72}d(J_{14},G)+\frac{4}{72}d(J_{20},G)+\frac{4}{72}d(J_{22},G)\]
\[t_{\inj}(B_{1}\cdot B_{8},G) = \frac{8}{72}d(J_{13},G)+\frac{4}{72}d(J_{19},G)+\frac{2}{72}d(J_{21},G)\]
\[t_{\inj}(B_{2}\cdot B_{2},G) = \frac{8}{72}d(J_{18},G)+\frac{4}{72}d(J_{21},G)+\frac{4}{72}d(J_{22},G)\]
\[t_{\inj}(B_{2}\cdot B_{3},G) = \frac{2}{72}d(J_{17},G)+\frac{2}{72}d(J_{20},G)+\frac{2}{72}d(J_{21},G)+\frac{4}{72}d(J_{22},G)\]
\[t_{\inj}(B_{2}\cdot B_{4},G) = \frac{2}{72}d(J_{16},G)+\frac{2}{72}d(J_{17},G)+\frac{2}{72}d(J_{19},G)+\frac{2}{72}d(J_{20},G)\]
\[t_{\inj}(B_{2}\cdot B_{5},G) = \frac{2}{72}d(J_{17},G)+\frac{2}{72}d(J_{20},G)+\frac{2}{72}d(J_{21},G)+\frac{4}{72}d(J_{22},G)\]
\[t_{\inj}(B_{2}\cdot B_{6},G) = \frac{2}{72}d(J_{16},G)+\frac{2}{72}d(J_{17},G)+\frac{2}{72}d(J_{19},G)+\frac{2}{72}d(J_{20},G)\]
\[t_{\inj}(B_{2}\cdot B_{7},G) = \frac{2}{72}d(J_{9},G)+\frac{4}{72}d(J_{11},G)+\frac{12}{72}d(J_{12},G)\]
\[t_{\inj}(B_{2}\cdot B_{8},G) = \frac{2}{72}d(J_{8},G)+\frac{4}{72}d(J_{10},G)+\frac{2}{72}d(J_{11},G)\]
\[t_{\inj}(B_{3}\cdot B_{3},G) = \frac{2}{72}d(J_{17},G)+\frac{4}{72}d(J_{21},G)+\frac{4}{72}d(J_{24},G)\]
\[t_{\inj}(B_{3}\cdot B_{4},G) = \frac{2}{72}d(J_{16},G)+\frac{8}{72}d(J_{18},G)+\frac{2}{72}d(J_{19},G)+\frac{2}{72}d(J_{21},G)\]
\[t_{\inj}(B_{3}\cdot B_{5},G) = \frac{2}{72}d(J_{11},G)+\frac{12}{72}d(J_{12},G)+\frac{2}{72}d(J_{21},G)+\frac{4}{72}d(J_{22},G)\]
\[t_{\inj}(B_{3}\cdot B_{6},G) = \frac{2}{72}d(J_{10},G)+\frac{2}{72}d(J_{11},G)+\frac{2}{72}d(J_{19},G)+\frac{2}{72}d(J_{20},G)\]
\[t_{\inj}(B_{3}\cdot B_{7},G) = \frac{2}{72}d(J_{9},G)+\frac{2}{72}d(J_{11},G)+\frac{2}{72}d(J_{17},G)+\frac{2}{72}d(J_{20},G)\]
\[t_{\inj}(B_{3}\cdot B_{8},G) = \frac{2}{72}d(J_{8},G)+\frac{2}{72}d(J_{10},G)+\frac{2}{72}d(J_{16},G)+\frac{2}{72}d(J_{17},G)\]
\[t_{\inj}(B_{4}\cdot B_{4},G) = \frac{12}{72}d(J_{15},G)+\frac{4}{72}d(J_{16},G)+\frac{2}{72}d(J_{17},G)\]
\[t_{\inj}(B_{4}\cdot B_{5},G) = \frac{2}{72}d(J_{10},G)+\frac{2}{72}d(J_{11},G)+\frac{2}{72}d(J_{19},G)+\frac{2}{72}d(J_{20},G)\]
\[t_{\inj}(B_{4}\cdot B_{6},G) = \frac{2}{72}d(J_{8},G)+\frac{2}{72}d(J_{9},G)+\frac{8}{72}d(J_{13},G)+\frac{8}{72}d(J_{14},G)\]
\[t_{\inj}(B_{4}\cdot B_{7},G) = \frac{2}{72}d(J_{6},G)+\frac{8}{72}d(J_{7},G)+\frac{2}{72}d(J_{10},G)+\frac{2}{72}d(J_{11},G)\]
\[t_{\inj}(B_{4}\cdot B_{8},G) = \frac{4}{72}d(J_{5},G)+\frac{2}{72}d(J_{6},G)+\frac{2}{72}d(J_{8},G)+\frac{2}{72}d(J_{9},G)\]
\[t_{\inj}(B_{5}\cdot B_{5},G) = \frac{2}{72}d(J_{17},G)+\frac{4}{72}d(J_{21},G)+\frac{4}{72}d(J_{24},G)\]
\[t_{\inj}(B_{5}\cdot B_{6},G) = \frac{2}{72}d(J_{16},G)+\frac{8}{72}d(J_{18},G)+\frac{2}{72}d(J_{19},G)+\frac{2}{72}d(J_{21},G)\]
\[t_{\inj}(B_{5}\cdot B_{7},G) = \frac{2}{72}d(J_{9},G)+\frac{2}{72}d(J_{11},G)+\frac{2}{72}d(J_{17},G)+\frac{2}{72}d(J_{20},G)\]
\[t_{\inj}(B_{5}\cdot B_{8},G) = \frac{2}{72}d(J_{8},G)+\frac{2}{72}d(J_{10},G)+\frac{2}{72}d(J_{16},G)+\frac{2}{72}d(J_{17},G)\]
\[t_{\inj}(B_{6}\cdot B_{6},G) = \frac{12}{72}d(J_{15},G)+\frac{4}{72}d(J_{16},G)+\frac{2}{72}d(J_{17},G)\]
\[t_{\inj}(B_{6}\cdot B_{7},G) = \frac{2}{72}d(J_{6},G)+\frac{8}{72}d(J_{7},G)+\frac{2}{72}d(J_{10},G)+\frac{2}{72}d(J_{11},G)\]
\[t_{\inj}(B_{6}\cdot B_{8},G) = \frac{4}{72}d(J_{5},G)+\frac{2}{72}d(J_{6},G)+\frac{2}{72}d(J_{8},G)+\frac{2}{72}d(J_{9},G)\]
\[t_{\inj}(B_{7}\cdot B_{7},G) = \frac{12}{72}d(J_{4},G)+\frac{4}{72}d(J_{9},G)+\frac{8}{72}d(J_{14},G)\]
\[t_{\inj}(B_{7}\cdot B_{8},G) = \frac{4}{72}d(J_{3},G)+\frac{4}{72}d(J_{6},G)+\frac{2}{72}d(J_{9},G)\]
\[t_{\inj}(B_{8}\cdot B_{8},G) = \frac{8}{72}d(J_{2},G)+\frac{8}{72}d(J_{3},G)+\frac{12}{72}d(J_{4},G).\]
Also, note that $t_{\inj}(R_{i}\cdot R_{j},G) = t_{\inj}(R_{j}\cdot R_{i},G)$ and $t_{\inj}(B_{i}\cdot B_{j},G) = t_{\inj}(B_{j}\cdot B_{i},G)$ for all $1\leq i,j\leq 8$. Our next step is to put all of this together to compute the coefficient of $d(J_\ell,G)$ in \eqref{eq:flagged} for each $1\leq \ell\leq 26$. In what follows, the calculation of the coefficient of $d(J_\ell,G)$ starts with $J_\ell$ followed by a colon. For each $\ell$, we first express the coefficient in terms of the entries of $A$ and then substitute the values from the matrix into the expression and simplify to obtain the final numerical coefficient. 
\begin{itemize}
\item[$J_{1}$:] $\displaystyle A(1,1)
$\\$\displaystyle =\frac{2}{128} = (1/2)^6$
\item[$J_{2}$:] $\displaystyle\frac{16}{72} A(1,1)+\frac{8}{72} A(1,2)+\frac{8}{72} A(1,3)+\frac{8}{72} A(1,5)+\frac{8}{72} A(2,1)+\frac{8}{72} A(3,1)$\\$\displaystyle+\frac{8}{72} A(5,1)+\frac{8}{72} A(8,8)
$\\$\displaystyle = \frac{16}{72} \left(\frac{2}{128}\right)+\frac{8}{72} \left(\frac{-6}{128}\right)+\frac{8}{72} \left(\frac{-2}{128}\right)+\frac{8}{72} \left(\frac{1}{128}\right)+\frac{8}{72} \left(\frac{-6}{128}\right)+\frac{8}{72} \left(\frac{-2}{128}\right)$\\$\displaystyle+\frac{8}{72} \left(\frac{1}{128}\right)+\frac{8}{72} \left(\frac{28}{128}\right) = (1/2)^6$
\item[$J_{3}$:] $\displaystyle\frac{4}{72} A(1,1)+\frac{4}{72} A(1,3)+\frac{4}{72} A(1,5)+\frac{4}{72} A(1,7)+\frac{4}{72} A(2,2)+\frac{4}{72} A(2,3)$\\$\displaystyle+\frac{4}{72} A(2,5)+\frac{4}{72} A(3,1)+\frac{4}{72} A(3,2)+\frac{4}{72} A(3,5)+\frac{4}{72} A(5,1)+\frac{4}{72} A(5,2)$\\$\displaystyle+\frac{4}{72} A(5,3)+\frac{4}{72} A(7,1)+\frac{4}{72} A(7,8)+\frac{4}{72} A(8,7)+\frac{8}{72} A(8,8)
$\\$\displaystyle = \frac{4}{72} \left(\frac{2}{128}\right)+\frac{4}{72} \left(\frac{-2}{128}\right)+\frac{4}{72} \left(\frac{1}{128}\right)+\frac{4}{72} \left(\frac{5}{128}\right)+\frac{4}{72} \left(\frac{58}{128}\right)+\frac{4}{72} \left(\frac{-3}{128}\right)$\\$\displaystyle+\frac{4}{72} \left(\frac{-6}{128}\right)+\frac{4}{72} \left(\frac{-2}{128}\right)+\frac{4}{72} \left(\frac{-3}{128}\right)+\frac{4}{72} \left(\frac{-47}{128}\right)+\frac{4}{72} \left(\frac{1}{128}\right)+\frac{4}{72} \left(\frac{-6}{128}\right)$\\$\displaystyle+\frac{4}{72} \left(\frac{-47}{128}\right)+\frac{4}{72} \left(\frac{5}{128}\right)+\frac{4}{72} \left(\frac{12}{128}\right)+\frac{4}{72} \left(\frac{12}{128}\right)+\frac{8}{72} \left(\frac{28}{128}\right) = (1/2)^6$
\item[$J_{4}$:] $\displaystyle\frac{1}{6}+\frac{12}{72} A(2,7)+\frac{12}{72} A(3,5)+\frac{12}{72} A(5,3)+\frac{12}{72} A(7,2)+\frac{12}{72} A(7,7)+\frac{12}{72} A(8,8)$\\$\displaystyle
$\\$\displaystyle =\frac{12}{72}+ \frac{12}{72} \left(\frac{-47}{128}\right)+\frac{12}{72} \left(\frac{-47}{128}\right)+\frac{12}{72} \left(\frac{-47}{128}\right)+\frac{12}{72} \left(\frac{-47}{128}\right)+\frac{12}{72} \left(\frac{44}{128}\right)+\frac{12}{72} \left(\frac{28}{128}\right)$\\$\displaystyle = (1/2)^6$
\item[$J_{5}$:] $\displaystyle\frac{8}{72} A(1,2)+\frac{4}{72} A(1,3)+\frac{4}{72} A(1,4)+\frac{4}{72} A(1,5)+\frac{4}{72} A(1,6)+\frac{8}{72} A(2,1)$\\$\displaystyle+\frac{4}{72} A(3,1)+\frac{4}{72} A(3,3)+\frac{4}{72} A(4,1)+\frac{4}{72} A(5,1)+\frac{4}{72} A(5,5)+\frac{4}{72} A(6,1)$\\$\displaystyle+\frac{4}{72} A(4,8)+\frac{4}{72} A(6,8)+\frac{4}{72} A(8,4)+\frac{4}{72} A(8,6)
$\\$\displaystyle = \frac{8}{72} \left(\frac{-6}{128}\right)+\frac{4}{72} \left(\frac{-2}{128}\right)+\frac{4}{72} \left(\frac{-3}{128}\right)+\frac{4}{72} \left(\frac{1}{128}\right)+\frac{4}{72} \left(\frac{-3}{128}\right)+\frac{8}{72} \left(\frac{-6}{128}\right)$\\$\displaystyle+\frac{4}{72} \left(\frac{-2}{128}\right)+\frac{4}{72} \left(\frac{56}{128}\right)+\frac{4}{72} \left(\frac{-3}{128}\right)+\frac{4}{72} \left(\frac{1}{128}\right)+\frac{4}{72} \left(\frac{56}{128}\right)+\frac{4}{72} \left(\frac{-3}{128}\right)$\\$\displaystyle+\frac{4}{72} \left(\frac{-9}{128}\right)+\frac{4}{72} \left(\frac{-10}{128}\right)+\frac{4}{72} \left(\frac{-9}{128}\right)+\frac{4}{72} \left(\frac{-10}{128}\right) = (1/2)^6$
\item[$J_{6}$:] $\displaystyle\frac{2}{72} A(1,3)+\frac{2}{72} A(1,5)+\frac{4}{72} A(1,7)+\frac{2}{72} A(2,3)+\frac{2}{72} A(2,4)+\frac{2}{72} A(2,5)$\\$\displaystyle+\frac{2}{72} A(2,6)+\frac{2}{72} A(3,1)+\frac{2}{72} A(3,2)+\frac{2}{72} A(3,6)+\frac{2}{72} A(3,7)+\frac{2}{72} A(4,2)$\\$\displaystyle+\frac{2}{72} A(4,5)+\frac{2}{72} A(5,1)+\frac{2}{72} A(5,2)+\frac{2}{72} A(5,4)+\frac{2}{72} A(5,7)+\frac{2}{72} A(6,2)$\\$\displaystyle+\frac{2}{72} A(6,3)+\frac{4}{72} A(7,1)+\frac{2}{72} A(7,3)+\frac{2}{72} A(7,5)+\frac{2}{72} A(4,7)+\frac{2}{72} A(4,8)$\\$\displaystyle+\frac{2}{72} A(6,7)+\frac{2}{72} A(6,8)+\frac{2}{72} A(7,4)+\frac{2}{72} A(7,6)+\frac{4}{72} A(7,8)+\frac{2}{72} A(8,4)$\\$\displaystyle+\frac{2}{72} A(8,6)+\frac{4}{72} A(8,7)
$\\$\displaystyle = \frac{2}{72} \left(\frac{-2}{128}\right)+\frac{2}{72} \left(\frac{1}{128}\right)+\frac{4}{72} \left(\frac{5}{128}\right)+\frac{2}{72} \left(\frac{-3}{128}\right)+\frac{2}{72} \left(\frac{12}{128}\right)+\frac{2}{72} \left(\frac{-6}{128}\right)$\\$\displaystyle+\frac{2}{72} \left(\frac{12}{128}\right)+\frac{2}{72} \left(\frac{-2}{128}\right)+\frac{2}{72} \left(\frac{-3}{128}\right)+\frac{2}{72} \left(\frac{11}{128}\right)+\frac{2}{72} \left(\frac{4}{128}\right)+\frac{2}{72} \left(\frac{12}{128}\right)$\\$\displaystyle+\frac{2}{72} \left(\frac{10}{128}\right)+\frac{2}{72} \left(\frac{1}{128}\right)+\frac{2}{72} \left(\frac{-6}{128}\right)+\frac{2}{72} \left(\frac{10}{128}\right)+\frac{2}{72} \left(\frac{2}{128}\right)+\frac{2}{72} \left(\frac{12}{128}\right)$\\$\displaystyle+\frac{2}{72} \left(\frac{11}{128}\right)+\frac{4}{72} \left(\frac{5}{128}\right)+\frac{2}{72} \left(\frac{4}{128}\right)+\frac{2}{72} \left(\frac{2}{128}\right)+\frac{2}{72} \left(\frac{-10}{128}\right)+\frac{2}{72} \left(\frac{-9}{128}\right)$\\$\displaystyle+\frac{2}{72} \left(\frac{-10}{128}\right)+\frac{2}{72} \left(\frac{-10}{128}\right)+\frac{2}{72} \left(\frac{-10}{128}\right)+\frac{2}{72} \left(\frac{-10}{128}\right)+\frac{4}{72} \left(\frac{12}{128}\right)+\frac{2}{72} \left(\frac{-9}{128}\right)$\\$\displaystyle+\frac{2}{72} \left(\frac{-10}{128}\right)+\frac{4}{72} \left(\frac{12}{128}\right) = (1/2)^6$
\item[$J_{7}$:] $\displaystyle\frac{8}{72} A(1,7)+\frac{8}{72} A(4,6)+\frac{8}{72} A(6,4)+\frac{8}{72} A(7,1)+\frac{8}{72} A(7,7)+\frac{8}{72} A(4,7)$\\$\displaystyle+\frac{8}{72} A(6,7)+\frac{8}{72} A(7,4)+\frac{8}{72} A(7,6)
$\\$\displaystyle = \frac{8}{72} \left(\frac{5}{128}\right)+\frac{8}{72} \left(\frac{2}{128}\right)+\frac{8}{72} \left(\frac{2}{128}\right)+\frac{8}{72} \left(\frac{5}{128}\right)+\frac{8}{72} \left(\frac{44}{128}\right)+\frac{8}{72} \left(\frac{-10}{128}\right)$\\$\displaystyle+\frac{8}{72} \left(\frac{-10}{128}\right)+\frac{8}{72} \left(\frac{-10}{128}\right)+\frac{8}{72} \left(\frac{-10}{128}\right) = (1/2)^6$
\item[$J_{8}$:] $\displaystyle\frac{2}{72} A(1,2)+\frac{2}{72} A(1,4)+\frac{2}{72} A(1,6)+\frac{2}{72} A(1,8)+\frac{2}{72} A(2,1)+\frac{4}{72} A(2,2)$\\$\displaystyle+\frac{2}{72} A(2,3)+\frac{2}{72} A(2,5)+\frac{2}{72} A(3,2)+\frac{4}{72} A(3,3)+\frac{2}{72} A(3,4)+\frac{2}{72} A(3,5)$\\$\displaystyle+\frac{2}{72} A(4,1)+\frac{2}{72} A(4,3)+\frac{2}{72} A(5,2)+\frac{2}{72} A(5,3)+\frac{4}{72} A(5,5)+\frac{2}{72} A(5,6)$\\$\displaystyle+\frac{2}{72} A(6,1)+\frac{2}{72} A(6,5)+\frac{2}{72} A(8,1)+\frac{2}{72} A(2,8)+\frac{2}{72} A(3,8)+\frac{2}{72} A(4,6)$\\$\displaystyle+\frac{2}{72} A(4,8)+\frac{2}{72} A(5,8)+\frac{2}{72} A(6,4)+\frac{2}{72} A(6,8)+\frac{2}{72} A(8,2)+\frac{2}{72} A(8,3)$\\$\displaystyle+\frac{2}{72} A(8,4)+\frac{2}{72} A(8,5)+\frac{2}{72} A(8,6)
$\\$\displaystyle = \frac{2}{72} \left(\frac{-6}{128}\right)+\frac{2}{72} \left(\frac{-3}{128}\right)+\frac{2}{72} \left(\frac{-3}{128}\right)+\frac{2}{72} \left(\frac{6}{128}\right)+\frac{2}{72} \left(\frac{-6}{128}\right)+\frac{4}{72} \left(\frac{58}{128}\right)$\\$\displaystyle+\frac{2}{72} \left(\frac{-3}{128}\right)+\frac{2}{72} \left(\frac{-6}{128}\right)+\frac{2}{72} \left(\frac{-3}{128}\right)+\frac{4}{72} \left(\frac{56}{128}\right)+\frac{2}{72} \left(\frac{-14}{128}\right)+\frac{2}{72} \left(\frac{-47}{128}\right)$\\$\displaystyle+\frac{2}{72} \left(\frac{-3}{128}\right)+\frac{2}{72} \left(\frac{-14}{128}\right)+\frac{2}{72} \left(\frac{-6}{128}\right)+\frac{2}{72} \left(\frac{-47}{128}\right)+\frac{4}{72} \left(\frac{56}{128}\right)+\frac{2}{72} \left(\frac{-14}{128}\right)$\\$\displaystyle+\frac{2}{72} \left(\frac{-3}{128}\right)+\frac{2}{72} \left(\frac{-14}{128}\right)+\frac{2}{72} \left(\frac{6}{128}\right)+\frac{2}{72} \left(\frac{-20}{128}\right)+\frac{2}{72} \left(\frac{-5}{128}\right)+\frac{2}{72} \left(\frac{2}{128}\right)$\\$\displaystyle+\frac{2}{72} \left(\frac{-9}{128}\right)+\frac{2}{72} \left(\frac{-2}{128}\right)+\frac{2}{72} \left(\frac{2}{128}\right)+\frac{2}{72} \left(\frac{-10}{128}\right)+\frac{2}{72} \left(\frac{-20}{128}\right)+\frac{2}{72} \left(\frac{-5}{128}\right)$\\$\displaystyle+\frac{2}{72} \left(\frac{-9}{128}\right)+\frac{2}{72} \left(\frac{-2}{128}\right)+\frac{2}{72} \left(\frac{-10}{128}\right) = (1/2)^6$
\item[$J_{9}$:] $\displaystyle\frac{1}{12}+\frac{4}{72} A(2,7)+\frac{2}{72} A(2,8)+\frac{2}{72} A(3,5)+\frac{2}{72} A(3,6)+\frac{2}{72} A(3,7)+\frac{2}{72} A(4,5)$\\$\displaystyle+\frac{2}{72} A(4,7)+\frac{2}{72} A(5,3)+\frac{2}{72} A(5,4)+\frac{2}{72} A(5,7)+\frac{2}{72} A(6,3)+\frac{2}{72} A(6,7)$\\$\displaystyle+\frac{4}{72} A(7,2)+\frac{2}{72} A(7,3)+\frac{2}{72} A(7,4)+\frac{2}{72} A(7,5)+\frac{2}{72} A(7,6)+\frac{2}{72} A(8,2)$\\$\displaystyle+\frac{2}{72} A(2,7)+\frac{2}{72} A(3,7)+\frac{2}{72} A(4,6)+\frac{2}{72} A(4,8)+\frac{2}{72} A(5,7)+\frac{2}{72} A(6,4)$\\$\displaystyle+\frac{2}{72} A(6,8)+\frac{2}{72} A(7,2)+\frac{2}{72} A(7,3)+\frac{2}{72} A(7,5)+\frac{4}{72} A(7,7)+\frac{2}{72} A(7,8)$\\$\displaystyle+\frac{2}{72} A(8,4)+\frac{2}{72} A(8,6)+\frac{2}{72} A(8,7)
$\\$\displaystyle = \frac{6}{72}+\frac{4}{72} \left(\frac{-47}{128}\right)+\frac{2}{72} \left(\frac{-20}{128}\right)+\frac{2}{72} \left(\frac{-47}{128}\right)+\frac{2}{72} \left(\frac{11}{128}\right)+\frac{2}{72} \left(\frac{4}{128}\right)+\frac{2}{72} \left(\frac{10}{128}\right)$\\$\displaystyle+\frac{2}{72} \left(\frac{-10}{128}\right)+\frac{2}{72} \left(\frac{-47}{128}\right)+\frac{2}{72} \left(\frac{10}{128}\right)+\frac{2}{72} \left(\frac{2}{128}\right)+\frac{2}{72} \left(\frac{11}{128}\right)+\frac{2}{72} \left(\frac{-10}{128}\right)$\\$\displaystyle+\frac{4}{72} \left(\frac{-47}{128}\right)+\frac{2}{72} \left(\frac{4}{128}\right)+\frac{2}{72} \left(\frac{-10}{128}\right)+\frac{2}{72} \left(\frac{2}{128}\right)+\frac{2}{72} \left(\frac{-10}{128}\right)+\frac{2}{72} \left(\frac{-20}{128}\right)$\\$\displaystyle+\frac{2}{72} \left(\frac{-47}{128}\right)+\frac{2}{72} \left(\frac{4}{128}\right)+\frac{2}{72} \left(\frac{2}{128}\right)+\frac{2}{72} \left(\frac{-9}{128}\right)+\frac{2}{72} \left(\frac{2}{128}\right)+\frac{2}{72} \left(\frac{2}{128}\right)$\\$\displaystyle+\frac{2}{72} \left(\frac{-10}{128}\right)+\frac{2}{72} \left(\frac{-47}{128}\right)+\frac{2}{72} \left(\frac{4}{128}\right)+\frac{2}{72} \left(\frac{2}{128}\right)+\frac{4}{72} \left(\frac{44}{128}\right)+\frac{2}{72} \left(\frac{12}{128}\right)$\\$\displaystyle+\frac{2}{72} \left(\frac{-9}{128}\right)+\frac{2}{72} \left(\frac{-10}{128}\right)+\frac{2}{72} \left(\frac{12}{128}\right) = (1/2)^6$
\item[$J_{10}$:] $\displaystyle\frac{2}{72} A(2,3)+\frac{2}{72} A(2,4)+\frac{2}{72} A(2,5)+\frac{2}{72} A(2,6)+\frac{2}{72} A(3,2)+\frac{2}{72} A(3,3)$\\$\displaystyle+\frac{2}{72} A(3,7)+\frac{2}{72} A(3,8)+\frac{2}{72} A(4,2)+\frac{2}{72} A(4,4)+\frac{2}{72} A(5,2)+\frac{2}{72} A(5,5)$\\$\displaystyle+\frac{2}{72} A(5,7)+\frac{2}{72} A(5,8)+\frac{2}{72} A(6,2)+\frac{2}{72} A(6,6)+\frac{2}{72} A(7,3)+\frac{2}{72} A(7,5)$\\$\displaystyle+\frac{2}{72} A(8,3)+\frac{2}{72} A(8,5)+\frac{4}{72} A(2,8)+\frac{2}{72} A(3,6)+\frac{2}{72} A(3,8)+\frac{2}{72} A(4,5)$\\$\displaystyle+\frac{2}{72} A(4,7)+\frac{2}{72} A(5,4)+\frac{2}{72} A(5,8)+\frac{2}{72} A(6,3)+\frac{2}{72} A(6,7)+\frac{2}{72} A(7,4)$\\$\displaystyle+\frac{2}{72} A(7,6)+\frac{4}{72} A(8,2)+\frac{2}{72} A(8,3)+\frac{2}{72} A(8,5)
$\\$\displaystyle = \frac{2}{72} \left(\frac{-3}{128}\right)+\frac{2}{72} \left(\frac{12}{128}\right)+\frac{2}{72} \left(\frac{-6}{128}\right)+\frac{2}{72} \left(\frac{12}{128}\right)+\frac{2}{72} \left(\frac{-3}{128}\right)+\frac{2}{72} \left(\frac{56}{128}\right)$\\$\displaystyle+\frac{2}{72} \left(\frac{4}{128}\right)+\frac{2}{72} \left(\frac{-5}{128}\right)+\frac{2}{72} \left(\frac{12}{128}\right)+\frac{2}{72} \left(\frac{12}{128}\right)+\frac{2}{72} \left(\frac{-6}{128}\right)+\frac{2}{72} \left(\frac{56}{128}\right)$\\$\displaystyle+\frac{2}{72} \left(\frac{2}{128}\right)+\frac{2}{72} \left(\frac{-2}{128}\right)+\frac{2}{72} \left(\frac{12}{128}\right)+\frac{2}{72} \left(\frac{12}{128}\right)+\frac{2}{72} \left(\frac{4}{128}\right)+\frac{2}{72} \left(\frac{2}{128}\right)$\\$\displaystyle+\frac{2}{72} \left(\frac{-5}{128}\right)+\frac{2}{72} \left(\frac{-2}{128}\right)+\frac{4}{72} \left(\frac{-20}{128}\right)+\frac{2}{72} \left(\frac{11}{128}\right)+\frac{2}{72} \left(\frac{-5}{128}\right)+\frac{2}{72} \left(\frac{10}{128}\right)$\\$\displaystyle+\frac{2}{72} \left(\frac{-10}{128}\right)+\frac{2}{72} \left(\frac{10}{128}\right)+\frac{2}{72} \left(\frac{-2}{128}\right)+\frac{2}{72} \left(\frac{11}{128}\right)+\frac{2}{72} \left(\frac{-10}{128}\right)+\frac{2}{72} \left(\frac{-10}{128}\right)$\\$\displaystyle+\frac{2}{72} \left(\frac{-10}{128}\right)+\frac{4}{72} \left(\frac{-20}{128}\right)+\frac{2}{72} \left(\frac{-5}{128}\right)+\frac{2}{72} \left(\frac{-2}{128}\right) = (1/2)^6$
\item[$J_{11}$:] $\displaystyle\frac{1}{12}+\frac{2}{72} A(2,7)+\frac{2}{72} A(3,7)+\frac{2}{72} A(4,6)+\frac{2}{72} A(4,8)+\frac{2}{72} A(5,7)+\frac{2}{72} A(6,4)$\\$\displaystyle+\frac{2}{72} A(6,8)+\frac{2}{72} A(7,2)+\frac{2}{72} A(7,3)+\frac{2}{72} A(7,5)+\frac{4}{72} A(7,7)+\frac{2}{72} A(7,8)$\\$\displaystyle+\frac{2}{72} A(8,4)+\frac{2}{72} A(8,6)+\frac{2}{72} A(8,7)+\frac{4}{72} A(2,7)+\frac{2}{72} A(2,8)+\frac{2}{72} A(3,5)$\\$\displaystyle+\frac{2}{72} A(3,6)+\frac{2}{72} A(3,7)+\frac{2}{72} A(4,5)+\frac{2}{72} A(4,7)+\frac{2}{72} A(5,3)+\frac{2}{72} A(5,4)$\\$\displaystyle+\frac{2}{72} A(5,7)+\frac{2}{72} A(6,3)+\frac{2}{72} A(6,7)+\frac{4}{72} A(7,2)+\frac{2}{72} A(7,3)+\frac{2}{72} A(7,4)$\\$\displaystyle+\frac{2}{72} A(7,5)+\frac{2}{72} A(7,6)+\frac{2}{72} A(8,2)
$\\$\displaystyle = \frac{6}{72}+\frac{2}{72} \left(\frac{-47}{128}\right)+\frac{2}{72} \left(\frac{4}{128}\right)+\frac{2}{72} \left(\frac{2}{128}\right)+\frac{2}{72} \left(\frac{-9}{128}\right)+\frac{2}{72} \left(\frac{2}{128}\right)+\frac{2}{72} \left(\frac{2}{128}\right)$\\$\displaystyle+\frac{2}{72} \left(\frac{-10}{128}\right)+\frac{2}{72} \left(\frac{-47}{128}\right)+\frac{2}{72} \left(\frac{4}{128}\right)+\frac{2}{72} \left(\frac{2}{128}\right)+\frac{4}{72} \left(\frac{44}{128}\right)+\frac{2}{72} \left(\frac{12}{128}\right)$\\$\displaystyle+\frac{2}{72} \left(\frac{-9}{128}\right)+\frac{2}{72} \left(\frac{-10}{128}\right)+\frac{2}{72} \left(\frac{12}{128}\right)+\frac{4}{72} \left(\frac{-47}{128}\right)+\frac{2}{72} \left(\frac{-20}{128}\right)+\frac{2}{72} \left(\frac{-47}{128}\right)$\\$\displaystyle+\frac{2}{72} \left(\frac{11}{128}\right)+\frac{2}{72} \left(\frac{4}{128}\right)+\frac{2}{72} \left(\frac{10}{128}\right)+\frac{2}{72} \left(\frac{-10}{128}\right)+\frac{2}{72} \left(\frac{-47}{128}\right)+\frac{2}{72} \left(\frac{10}{128}\right)$\\$\displaystyle+\frac{2}{72} \left(\frac{2}{128}\right)+\frac{2}{72} \left(\frac{11}{128}\right)+\frac{2}{72} \left(\frac{-10}{128}\right)+\frac{4}{72} \left(\frac{-47}{128}\right)+\frac{2}{72} \left(\frac{4}{128}\right)+\frac{2}{72} \left(\frac{-10}{128}\right)$\\$\displaystyle+\frac{2}{72} \left(\frac{2}{128}\right)+\frac{2}{72} \left(\frac{-10}{128}\right)+\frac{2}{72} \left(\frac{-20}{128}\right) = (1/2)^6$
\item[$J_{12}$:] $\displaystyle\frac{1}{6}+\frac{12}{72} A(7,7)+\frac{12}{72} A(8,8)+\frac{12}{72} A(2,7)+\frac{12}{72} A(3,5)+\frac{12}{72} A(5,3)+\frac{12}{72} A(7,2)$\\$\displaystyle
$\\$\displaystyle = \frac{12}{72}+\frac{12}{72} \left(\frac{44}{128}\right)+\frac{12}{72} \left(\frac{28}{128}\right)+\frac{12}{72} \left(\frac{-47}{128}\right)+\frac{12}{72} \left(\frac{-47}{128}\right)+\frac{12}{72} \left(\frac{-47}{128}\right)+\frac{12}{72} \left(\frac{-47}{128}\right)$\\$\displaystyle = (1/2)^6$
\item[$J_{13}$:] $\displaystyle\frac{8}{72} A(2,2)+\frac{8}{72} A(3,4)+\frac{8}{72} A(4,3)+\frac{8}{72} A(5,6)+\frac{8}{72} A(6,5)+\frac{8}{72} A(1,8)$\\$\displaystyle+\frac{8}{72} A(4,6)+\frac{8}{72} A(6,4)+\frac{8}{72} A(8,1)
$\\$\displaystyle = \frac{8}{72} \left(\frac{58}{128}\right)+\frac{8}{72} \left(\frac{-14}{128}\right)+\frac{8}{72} \left(\frac{-14}{128}\right)+\frac{8}{72} \left(\frac{-14}{128}\right)+\frac{8}{72} \left(\frac{-14}{128}\right)+\frac{8}{72} \left(\frac{6}{128}\right)$\\$\displaystyle+\frac{8}{72} \left(\frac{2}{128}\right)+\frac{8}{72} \left(\frac{2}{128}\right)+\frac{8}{72} \left(\frac{6}{128}\right) = (1/2)^6$
\item[$J_{14}$:] $\displaystyle\frac{8}{72} A(4,7)+\frac{8}{72} A(6,7)+\frac{8}{72} A(7,4)+\frac{8}{72} A(7,6)+\frac{8}{72} A(1,7)+\frac{8}{72} A(4,6)$\\$\displaystyle+\frac{8}{72} A(6,4)+\frac{8}{72} A(7,1)+\frac{8}{72} A(7,7)
$\\$\displaystyle = \frac{8}{72} \left(\frac{-10}{128}\right)+\frac{8}{72} \left(\frac{-10}{128}\right)+\frac{8}{72} \left(\frac{-10}{128}\right)+\frac{8}{72} \left(\frac{-10}{128}\right)+\frac{8}{72} \left(\frac{5}{128}\right)+\frac{8}{72} \left(\frac{2}{128}\right)$\\$\displaystyle+\frac{8}{72} \left(\frac{2}{128}\right)+\frac{8}{72} \left(\frac{5}{128}\right)+\frac{8}{72} \left(\frac{44}{128}\right) = (1/2)^6$
\item[$J_{15}$:] $\displaystyle\frac{12}{72} A(1,4)+\frac{12}{72} A(1,6)+\frac{12}{72} A(4,1)+\frac{12}{72} A(6,1)+\frac{12}{72} A(4,4)+\frac{12}{72} A(6,6)$\\$\displaystyle
$\\$\displaystyle = \frac{12}{72} \left(\frac{-3}{128}\right)+\frac{12}{72} \left(\frac{-3}{128}\right)+\frac{12}{72} \left(\frac{-3}{128}\right)+\frac{12}{72} \left(\frac{-3}{128}\right)+\frac{12}{72} \left(\frac{12}{128}\right)+\frac{12}{72} \left(\frac{12}{128}\right)$\\$\displaystyle = (1/2)^6$
\item[$J_{16}$:] $\displaystyle\frac{2}{72} A(1,4)+\frac{2}{72} A(1,6)+\frac{4}{72} A(1,8)+\frac{2}{72} A(2,4)+\frac{2}{72} A(2,6)+\frac{2}{72} A(3,4)$\\$\displaystyle+\frac{2}{72} A(3,6)+\frac{2}{72} A(4,1)+\frac{2}{72} A(4,2)+\frac{2}{72} A(4,3)+\frac{2}{72} A(4,5)+\frac{2}{72} A(5,4)$\\$\displaystyle+\frac{2}{72} A(5,6)+\frac{2}{72} A(6,1)+\frac{2}{72} A(6,2)+\frac{2}{72} A(6,3)+\frac{2}{72} A(6,5)+\frac{4}{72} A(8,1)$\\$\displaystyle+\frac{2}{72} A(2,4)+\frac{2}{72} A(2,6)+\frac{2}{72} A(3,4)+\frac{2}{72} A(3,8)+\frac{2}{72} A(4,2)+\frac{2}{72} A(4,3)$\\$\displaystyle+\frac{4}{72} A(4,4)+\frac{2}{72} A(5,6)+\frac{2}{72} A(5,8)+\frac{2}{72} A(6,2)+\frac{2}{72} A(6,5)+\frac{4}{72} A(6,6)$\\$\displaystyle+\frac{2}{72} A(8,3)+\frac{2}{72} A(8,5)
$\\$\displaystyle = \frac{2}{72} \left(\frac{-3}{128}\right)+\frac{2}{72} \left(\frac{-3}{128}\right)+\frac{4}{72} \left(\frac{6}{128}\right)+\frac{2}{72} \left(\frac{12}{128}\right)+\frac{2}{72} \left(\frac{12}{128}\right)+\frac{2}{72} \left(\frac{-14}{128}\right)$\\$\displaystyle+\frac{2}{72} \left(\frac{11}{128}\right)+\frac{2}{72} \left(\frac{-3}{128}\right)+\frac{2}{72} \left(\frac{12}{128}\right)+\frac{2}{72} \left(\frac{-14}{128}\right)+\frac{2}{72} \left(\frac{10}{128}\right)+\frac{2}{72} \left(\frac{10}{128}\right)$\\$\displaystyle+\frac{2}{72} \left(\frac{-14}{128}\right)+\frac{2}{72} \left(\frac{-3}{128}\right)+\frac{2}{72} \left(\frac{12}{128}\right)+\frac{2}{72} \left(\frac{11}{128}\right)+\frac{2}{72} \left(\frac{-14}{128}\right)+\frac{4}{72} \left(\frac{6}{128}\right)$\\$\displaystyle+\frac{2}{72} \left(\frac{12}{128}\right)+\frac{2}{72} \left(\frac{12}{128}\right)+\frac{2}{72} \left(\frac{-14}{128}\right)+\frac{2}{72} \left(\frac{-5}{128}\right)+\frac{2}{72} \left(\frac{12}{128}\right)+\frac{2}{72} \left(\frac{-14}{128}\right)$\\$\displaystyle+\frac{4}{72} \left(\frac{12}{128}\right)+\frac{2}{72} \left(\frac{-14}{128}\right)+\frac{2}{72} \left(\frac{-2}{128}\right)+\frac{2}{72} \left(\frac{12}{128}\right)+\frac{2}{72} \left(\frac{-14}{128}\right)+\frac{4}{72} \left(\frac{12}{128}\right)$\\$\displaystyle+\frac{2}{72} \left(\frac{-5}{128}\right)+\frac{2}{72} \left(\frac{-2}{128}\right) = (1/2)^6$
\item[$J_{17}$:] $\displaystyle\frac{4}{72} A(2,8)+\frac{2}{72} A(3,6)+\frac{2}{72} A(3,8)+\frac{2}{72} A(4,5)+\frac{2}{72} A(4,7)+\frac{2}{72} A(5,4)$\\$\displaystyle+\frac{2}{72} A(5,8)+\frac{2}{72} A(6,3)+\frac{2}{72} A(6,7)+\frac{2}{72} A(7,4)+\frac{2}{72} A(7,6)+\frac{4}{72} A(8,2)$\\$\displaystyle+\frac{2}{72} A(8,3)+\frac{2}{72} A(8,5)+\frac{2}{72} A(2,3)+\frac{2}{72} A(2,4)+\frac{2}{72} A(2,5)+\frac{2}{72} A(2,6)$\\$\displaystyle+\frac{2}{72} A(3,2)+\frac{2}{72} A(3,3)+\frac{2}{72} A(3,7)+\frac{2}{72} A(3,8)+\frac{2}{72} A(4,2)+\frac{2}{72} A(4,4)$\\$\displaystyle+\frac{2}{72} A(5,2)+\frac{2}{72} A(5,5)+\frac{2}{72} A(5,7)+\frac{2}{72} A(5,8)+\frac{2}{72} A(6,2)+\frac{2}{72} A(6,6)$\\$\displaystyle+\frac{2}{72} A(7,3)+\frac{2}{72} A(7,5)+\frac{2}{72} A(8,3)+\frac{2}{72} A(8,5)
$\\$\displaystyle = \frac{4}{72} \left(\frac{-20}{128}\right)+\frac{2}{72} \left(\frac{11}{128}\right)+\frac{2}{72} \left(\frac{-5}{128}\right)+\frac{2}{72} \left(\frac{10}{128}\right)+\frac{2}{72} \left(\frac{-10}{128}\right)+\frac{2}{72} \left(\frac{10}{128}\right)$\\$\displaystyle+\frac{2}{72} \left(\frac{-2}{128}\right)+\frac{2}{72} \left(\frac{11}{128}\right)+\frac{2}{72} \left(\frac{-10}{128}\right)+\frac{2}{72} \left(\frac{-10}{128}\right)+\frac{2}{72} \left(\frac{-10}{128}\right)+\frac{4}{72} \left(\frac{-20}{128}\right)$\\$\displaystyle+\frac{2}{72} \left(\frac{-5}{128}\right)+\frac{2}{72} \left(\frac{-2}{128}\right)+\frac{2}{72} \left(\frac{-3}{128}\right)+\frac{2}{72} \left(\frac{12}{128}\right)+\frac{2}{72} \left(\frac{-6}{128}\right)+\frac{2}{72} \left(\frac{12}{128}\right)$\\$\displaystyle+\frac{2}{72} \left(\frac{-3}{128}\right)+\frac{2}{72} \left(\frac{56}{128}\right)+\frac{2}{72} \left(\frac{4}{128}\right)+\frac{2}{72} \left(\frac{-5}{128}\right)+\frac{2}{72} \left(\frac{12}{128}\right)+\frac{2}{72} \left(\frac{12}{128}\right)$\\$\displaystyle+\frac{2}{72} \left(\frac{-6}{128}\right)+\frac{2}{72} \left(\frac{56}{128}\right)+\frac{2}{72} \left(\frac{2}{128}\right)+\frac{2}{72} \left(\frac{-2}{128}\right)+\frac{2}{72} \left(\frac{12}{128}\right)+\frac{2}{72} \left(\frac{12}{128}\right)$\\$\displaystyle+\frac{2}{72} \left(\frac{4}{128}\right)+\frac{2}{72} \left(\frac{2}{128}\right)+\frac{2}{72} \left(\frac{-5}{128}\right)+\frac{2}{72} \left(\frac{-2}{128}\right) = (1/2)^6$
\item[$J_{18}$:] $\displaystyle\frac{8}{72} A(1,8)+\frac{8}{72} A(4,6)+\frac{8}{72} A(6,4)+\frac{8}{72} A(8,1)+\frac{8}{72} A(2,2)+\frac{8}{72} A(3,4)$\\$\displaystyle+\frac{8}{72} A(4,3)+\frac{8}{72} A(5,6)+\frac{8}{72} A(6,5)
$\\$\displaystyle = \frac{8}{72} \left(\frac{6}{128}\right)+\frac{8}{72} \left(\frac{2}{128}\right)+\frac{8}{72} \left(\frac{2}{128}\right)+\frac{8}{72} \left(\frac{6}{128}\right)+\frac{8}{72} \left(\frac{58}{128}\right)+\frac{8}{72} \left(\frac{-14}{128}\right)$\\$\displaystyle+\frac{8}{72} \left(\frac{-14}{128}\right)+\frac{8}{72} \left(\frac{-14}{128}\right)+\frac{8}{72} \left(\frac{-14}{128}\right) = (1/2)^6$
\item[$J_{19}$:] $\displaystyle\frac{2}{72} A(2,4)+\frac{2}{72} A(2,6)+\frac{2}{72} A(3,4)+\frac{2}{72} A(3,8)+\frac{2}{72} A(4,2)+\frac{2}{72} A(4,3)$\\$\displaystyle+\frac{4}{72} A(4,4)+\frac{2}{72} A(5,6)+\frac{2}{72} A(5,8)+\frac{2}{72} A(6,2)+\frac{2}{72} A(6,5)+\frac{4}{72} A(6,6)$\\$\displaystyle+\frac{2}{72} A(8,3)+\frac{2}{72} A(8,5)+\frac{2}{72} A(1,4)+\frac{2}{72} A(1,6)+\frac{4}{72} A(1,8)+\frac{2}{72} A(2,4)$\\$\displaystyle+\frac{2}{72} A(2,6)+\frac{2}{72} A(3,4)+\frac{2}{72} A(3,6)+\frac{2}{72} A(4,1)+\frac{2}{72} A(4,2)+\frac{2}{72} A(4,3)$\\$\displaystyle+\frac{2}{72} A(4,5)+\frac{2}{72} A(5,4)+\frac{2}{72} A(5,6)+\frac{2}{72} A(6,1)+\frac{2}{72} A(6,2)+\frac{2}{72} A(6,3)$\\$\displaystyle+\frac{2}{72} A(6,5)+\frac{4}{72} A(8,1)
$\\$\displaystyle = \frac{2}{72} \left(\frac{12}{128}\right)+\frac{2}{72} \left(\frac{12}{128}\right)+\frac{2}{72} \left(\frac{-14}{128}\right)+\frac{2}{72} \left(\frac{-5}{128}\right)+\frac{2}{72} \left(\frac{12}{128}\right)+\frac{2}{72} \left(\frac{-14}{128}\right)$\\$\displaystyle+\frac{4}{72} \left(\frac{12}{128}\right)+\frac{2}{72} \left(\frac{-14}{128}\right)+\frac{2}{72} \left(\frac{-2}{128}\right)+\frac{2}{72} \left(\frac{12}{128}\right)+\frac{2}{72} \left(\frac{-14}{128}\right)+\frac{4}{72} \left(\frac{12}{128}\right)$\\$\displaystyle+\frac{2}{72} \left(\frac{-5}{128}\right)+\frac{2}{72} \left(\frac{-2}{128}\right)+\frac{2}{72} \left(\frac{-3}{128}\right)+\frac{2}{72} \left(\frac{-3}{128}\right)+\frac{4}{72} \left(\frac{6}{128}\right)+\frac{2}{72} \left(\frac{12}{128}\right)$\\$\displaystyle+\frac{2}{72} \left(\frac{12}{128}\right)+\frac{2}{72} \left(\frac{-14}{128}\right)+\frac{2}{72} \left(\frac{11}{128}\right)+\frac{2}{72} \left(\frac{-3}{128}\right)+\frac{2}{72} \left(\frac{12}{128}\right)+\frac{2}{72} \left(\frac{-14}{128}\right)$\\$\displaystyle+\frac{2}{72} \left(\frac{10}{128}\right)+\frac{2}{72} \left(\frac{10}{128}\right)+\frac{2}{72} \left(\frac{-14}{128}\right)+\frac{2}{72} \left(\frac{-3}{128}\right)+\frac{2}{72} \left(\frac{12}{128}\right)+\frac{2}{72} \left(\frac{11}{128}\right)$\\$\displaystyle+\frac{2}{72} \left(\frac{-14}{128}\right)+\frac{4}{72} \left(\frac{6}{128}\right) = (1/2)^6$
\item[$J_{20}$:] $\displaystyle\frac{2}{72} A(4,7)+\frac{2}{72} A(4,8)+\frac{2}{72} A(6,7)+\frac{2}{72} A(6,8)+\frac{2}{72} A(7,4)+\frac{2}{72} A(7,6)$\\$\displaystyle+\frac{4}{72} A(7,8)+\frac{2}{72} A(8,4)+\frac{2}{72} A(8,6)+\frac{4}{72} A(8,7)+\frac{2}{72} A(1,3)+\frac{2}{72} A(1,5)$\\$\displaystyle+\frac{4}{72} A(1,7)+\frac{2}{72} A(2,3)+\frac{2}{72} A(2,4)+\frac{2}{72} A(2,5)+\frac{2}{72} A(2,6)+\frac{2}{72} A(3,1)$\\$\displaystyle+\frac{2}{72} A(3,2)+\frac{2}{72} A(3,6)+\frac{2}{72} A(3,7)+\frac{2}{72} A(4,2)+\frac{2}{72} A(4,5)+\frac{2}{72} A(5,1)$\\$\displaystyle+\frac{2}{72} A(5,2)+\frac{2}{72} A(5,4)+\frac{2}{72} A(5,7)+\frac{2}{72} A(6,2)+\frac{2}{72} A(6,3)+\frac{4}{72} A(7,1)$\\$\displaystyle+\frac{2}{72} A(7,3)+\frac{2}{72} A(7,5)
$\\$\displaystyle = \frac{2}{72} \left(\frac{-10}{128}\right)+\frac{2}{72} \left(\frac{-9}{128}\right)+\frac{2}{72} \left(\frac{-10}{128}\right)+\frac{2}{72} \left(\frac{-10}{128}\right)+\frac{2}{72} \left(\frac{-10}{128}\right)+\frac{2}{72} \left(\frac{-10}{128}\right)$\\$\displaystyle+\frac{4}{72} \left(\frac{12}{128}\right)+\frac{2}{72} \left(\frac{-9}{128}\right)+\frac{2}{72} \left(\frac{-10}{128}\right)+\frac{4}{72} \left(\frac{12}{128}\right)+\frac{2}{72} \left(\frac{-2}{128}\right)+\frac{2}{72} \left(\frac{1}{128}\right)$\\$\displaystyle+\frac{4}{72} \left(\frac{5}{128}\right)+\frac{2}{72} \left(\frac{-3}{128}\right)+\frac{2}{72} \left(\frac{12}{128}\right)+\frac{2}{72} \left(\frac{-6}{128}\right)+\frac{2}{72} \left(\frac{12}{128}\right)+\frac{2}{72} \left(\frac{-2}{128}\right)$\\$\displaystyle+\frac{2}{72} \left(\frac{-3}{128}\right)+\frac{2}{72} \left(\frac{11}{128}\right)+\frac{2}{72} \left(\frac{4}{128}\right)+\frac{2}{72} \left(\frac{12}{128}\right)+\frac{2}{72} \left(\frac{10}{128}\right)+\frac{2}{72} \left(\frac{1}{128}\right)$\\$\displaystyle+\frac{2}{72} \left(\frac{-6}{128}\right)+\frac{2}{72} \left(\frac{10}{128}\right)+\frac{2}{72} \left(\frac{2}{128}\right)+\frac{2}{72} \left(\frac{12}{128}\right)+\frac{2}{72} \left(\frac{11}{128}\right)+\frac{4}{72} \left(\frac{5}{128}\right)$\\$\displaystyle+\frac{2}{72} \left(\frac{4}{128}\right)+\frac{2}{72} \left(\frac{2}{128}\right) = (1/2)^6$
\item[$J_{21}$:] $\displaystyle\frac{2}{72} A(2,8)+\frac{2}{72} A(3,8)+\frac{2}{72} A(4,6)+\frac{2}{72} A(4,8)+\frac{2}{72} A(5,8)+\frac{2}{72} A(6,4)$\\$\displaystyle+\frac{2}{72} A(6,8)+\frac{2}{72} A(8,2)+\frac{2}{72} A(8,3)+\frac{2}{72} A(8,4)+\frac{2}{72} A(8,5)+\frac{2}{72} A(8,6)$\\$\displaystyle+\frac{2}{72} A(1,2)+\frac{2}{72} A(1,4)+\frac{2}{72} A(1,6)+\frac{2}{72} A(1,8)+\frac{2}{72} A(2,1)+\frac{4}{72} A(2,2)$\\$\displaystyle+\frac{2}{72} A(2,3)+\frac{2}{72} A(2,5)+\frac{2}{72} A(3,2)+\frac{4}{72} A(3,3)+\frac{2}{72} A(3,4)+\frac{2}{72} A(3,5)$\\$\displaystyle+\frac{2}{72} A(4,1)+\frac{2}{72} A(4,3)+\frac{2}{72} A(5,2)+\frac{2}{72} A(5,3)+\frac{4}{72} A(5,5)+\frac{2}{72} A(5,6)$\\$\displaystyle+\frac{2}{72} A(6,1)+\frac{2}{72} A(6,5)+\frac{2}{72} A(8,1)
$\\$\displaystyle = \frac{2}{72} \left(\frac{-20}{128}\right)+\frac{2}{72} \left(\frac{-5}{128}\right)+\frac{2}{72} \left(\frac{2}{128}\right)+\frac{2}{72} \left(\frac{-9}{128}\right)+\frac{2}{72} \left(\frac{-2}{128}\right)+\frac{2}{72} \left(\frac{2}{128}\right)$\\$\displaystyle+\frac{2}{72} \left(\frac{-10}{128}\right)+\frac{2}{72} \left(\frac{-20}{128}\right)+\frac{2}{72} \left(\frac{-5}{128}\right)+\frac{2}{72} \left(\frac{-9}{128}\right)+\frac{2}{72} \left(\frac{-2}{128}\right)+\frac{2}{72} \left(\frac{-10}{128}\right)$\\$\displaystyle+\frac{2}{72} \left(\frac{-6}{128}\right)+\frac{2}{72} \left(\frac{-3}{128}\right)+\frac{2}{72} \left(\frac{-3}{128}\right)+\frac{2}{72} \left(\frac{6}{128}\right)+\frac{2}{72} \left(\frac{-6}{128}\right)+\frac{4}{72} \left(\frac{58}{128}\right)$\\$\displaystyle+\frac{2}{72} \left(\frac{-3}{128}\right)+\frac{2}{72} \left(\frac{-6}{128}\right)+\frac{2}{72} \left(\frac{-3}{128}\right)+\frac{4}{72} \left(\frac{56}{128}\right)+\frac{2}{72} \left(\frac{-14}{128}\right)+\frac{2}{72} \left(\frac{-47}{128}\right)$\\$\displaystyle+\frac{2}{72} \left(\frac{-3}{128}\right)+\frac{2}{72} \left(\frac{-14}{128}\right)+\frac{2}{72} \left(\frac{-6}{128}\right)+\frac{2}{72} \left(\frac{-47}{128}\right)+\frac{4}{72} \left(\frac{56}{128}\right)+\frac{2}{72} \left(\frac{-14}{128}\right)$\\$\displaystyle+\frac{2}{72} \left(\frac{-3}{128}\right)+\frac{2}{72} \left(\frac{-14}{128}\right)+\frac{2}{72} \left(\frac{6}{128}\right) = (1/2)^6$
\item[$J_{22}$:] $\displaystyle\frac{4}{72} A(7,8)+\frac{4}{72} A(8,7)+\frac{8}{72} A(8,8)+\frac{4}{72} A(1,1)+\frac{4}{72} A(1,3)+\frac{4}{72} A(1,5)$\\$\displaystyle+\frac{4}{72} A(1,7)+\frac{4}{72} A(2,2)+\frac{4}{72} A(2,3)+\frac{4}{72} A(2,5)+\frac{4}{72} A(3,1)+\frac{4}{72} A(3,2)$\\$\displaystyle+\frac{4}{72} A(3,5)+\frac{4}{72} A(5,1)+\frac{4}{72} A(5,2)+\frac{4}{72} A(5,3)+\frac{4}{72} A(7,1)
$\\$\displaystyle = \frac{4}{72} \left(\frac{12}{128}\right)+\frac{4}{72} \left(\frac{12}{128}\right)+\frac{8}{72} \left(\frac{28}{128}\right)+\frac{4}{72} \left(\frac{2}{128}\right)+\frac{4}{72} \left(\frac{-2}{128}\right)+\frac{4}{72} \left(\frac{1}{128}\right)$\\$\displaystyle+\frac{4}{72} \left(\frac{5}{128}\right)+\frac{4}{72} \left(\frac{58}{128}\right)+\frac{4}{72} \left(\frac{-3}{128}\right)+\frac{4}{72} \left(\frac{-6}{128}\right)+\frac{4}{72} \left(\frac{-2}{128}\right)+\frac{4}{72} \left(\frac{-3}{128}\right)$\\$\displaystyle+\frac{4}{72} \left(\frac{-47}{128}\right)+\frac{4}{72} \left(\frac{1}{128}\right)+\frac{4}{72} \left(\frac{-6}{128}\right)+\frac{4}{72} \left(\frac{-47}{128}\right)+\frac{4}{72} \left(\frac{5}{128}\right) = (1/2)^6$
\item[$J_{23}$:] $\displaystyle\frac{12}{72} A(4,4)+\frac{12}{72} A(6,6)+\frac{12}{72} A(1,4)+\frac{12}{72} A(1,6)+\frac{12}{72} A(4,1)+\frac{12}{72} A(6,1)$\\$\displaystyle
$\\$\displaystyle = \frac{12}{72} \left(\frac{12}{128}\right)+\frac{12}{72} \left(\frac{12}{128}\right)+\frac{12}{72} \left(\frac{-3}{128}\right)+\frac{12}{72} \left(\frac{-3}{128}\right)+\frac{12}{72} \left(\frac{-3}{128}\right)+\frac{12}{72} \left(\frac{-3}{128}\right)$\\$\displaystyle = (1/2)^6$
\item[$J_{24}$:] $\displaystyle\frac{4}{72} A(4,8)+\frac{4}{72} A(6,8)+\frac{4}{72} A(8,4)+\frac{4}{72} A(8,6)+\frac{8}{72} A(1,2)+\frac{4}{72} A(1,3)$\\$\displaystyle+\frac{4}{72} A(1,4)+\frac{4}{72} A(1,5)+\frac{4}{72} A(1,6)+\frac{8}{72} A(2,1)+\frac{4}{72} A(3,1)+\frac{4}{72} A(3,3)$\\$\displaystyle+\frac{4}{72} A(4,1)+\frac{4}{72} A(5,1)+\frac{4}{72} A(5,5)+\frac{4}{72} A(6,1)
$\\$\displaystyle = \frac{4}{72} \left(\frac{-9}{128}\right)+\frac{4}{72} \left(\frac{-10}{128}\right)+\frac{4}{72} \left(\frac{-9}{128}\right)+\frac{4}{72} \left(\frac{-10}{128}\right)+\frac{8}{72} \left(\frac{-6}{128}\right)+\frac{4}{72} \left(\frac{-2}{128}\right)$\\$\displaystyle+\frac{4}{72} \left(\frac{-3}{128}\right)+\frac{4}{72} \left(\frac{1}{128}\right)+\frac{4}{72} \left(\frac{-3}{128}\right)+\frac{8}{72} \left(\frac{-6}{128}\right)+\frac{4}{72} \left(\frac{-2}{128}\right)+\frac{4}{72} \left(\frac{56}{128}\right)$\\$\displaystyle+\frac{4}{72} \left(\frac{-3}{128}\right)+\frac{4}{72} \left(\frac{1}{128}\right)+\frac{4}{72} \left(\frac{56}{128}\right)+\frac{4}{72} \left(\frac{-3}{128}\right) = (1/2)^6$
\item[$J_{25}$:] $\displaystyle\frac{8}{72} A(8,8)+\frac{16}{72} A(1,1)+\frac{8}{72} A(1,2)+\frac{8}{72} A(1,3)+\frac{8}{72} A(1,5)+\frac{8}{72} A(2,1)$\\$\displaystyle+\frac{8}{72} A(3,1)+\frac{8}{72} A(5,1)
$\\$\displaystyle = \frac{8}{72} \left(\frac{28}{128}\right)+\frac{16}{72} \left(\frac{2}{128}\right)+\frac{8}{72} \left(\frac{-6}{128}\right)+\frac{8}{72} \left(\frac{-2}{128}\right)+\frac{8}{72} \left(\frac{1}{128}\right)+\frac{8}{72} \left(\frac{-6}{128}\right)$\\$\displaystyle+\frac{8}{72} \left(\frac{-2}{128}\right)+\frac{8}{72} \left(\frac{1}{128}\right) = (1/2)^6$
\item[$J_{26}$:] $A(1,1)$\\
$\displaystyle =  \frac{2}{128} = (1/2)^6$
\end{itemize}
So, we have that the expression
\[
\begin{aligned}
&\frac{1}{6}d(J_4,G) + \frac{1}{12}d(J_9,G)
 + \frac{1}{12}d(J_{11},G) + \frac{1}{6}d(J_{12},G)\\
&\quad + \sum_{i=1}^8\sum_{j=1}^8 A(i,j)
\left(
t_{\inj}(R_i\cdot R_j,G)
+t_{\inj}(B_i\cdot B_j,G)
\right)
\end{aligned}
\]
is equal to
\[\sum_{\ell=1}^{26}(1/2)^6d(J_{\ell},G)= (1/2)^6,\]
which completes the proof.
\end{proof}

\end{document}